\documentclass[11pt]{article}

\usepackage{fullwidth}
\usepackage{xypic}
\usepackage{xy}
\usepackage[top=1.5in, bottom=1.5in, left=1in, right=1in]{geometry}
\usepackage{amsgen}
\usepackage{amsmath}
\usepackage{amstext}
\usepackage{amsbsy}
\usepackage{amsopn}
\usepackage{amsfonts}
\usepackage{amssymb}
\usepackage{eepic}
\usepackage{graphicx}
\usepackage{epsf}
\usepackage{pstricks}
\usepackage{hyperref}

\xyoption{all}

\def\Box{\square}

\def\tra#1{\smash{\mathop{\mid\kern
-1pt\joinrel\relbar\joinrel\relbar}\limits^{*}_{#1}}}
\def\longtra#1{\smash{\mathop{\mid\kern
-1pt\joinrel\relbar\joinrel\relbar\joinrel\relbar}\limits^{*}_{#1}}}
\def\vlongtra#1{\smash{\mathop{\mid\kern
-1pt\joinrel\relbar\joinrel\relbar\joinrel\relbar\joinrel\relbar}\limits^{*}_{#1}}}
\def\vvlongtra#1{\smash{\mathop{\mid\kern
-1pt\joinrel\relbar\joinrel\relbar\joinrel\relbar\joinrel\relbar\joinrel\relbar}\limits^{*}_{#1}}}
\def\vvvlongtra#1{\smash{\mathop{\mid\kern
-1pt\joinrel\relbar\joinrel\relbar\joinrel\relbar\joinrel\relbar\joinrel\relbar\joinrel\relbar}\limits^{*}_{#1}}}
\def\etra#1{\smash{\mathop{\mid\kern
-1pt\joinrel\relbar\joinrel\relbar}\limits_{#1}}}

\def\A{{\cal{A}}}

\def\Rw{\Rightarrow}
\def\oo{\overline}

\def\M{{\bf M}}
\def\N{\mathbb{N}}
\def\U{{\bf U}}
\def\ab{{\bf Ab}}
\def\bu{{\bf B}}

\def\aut{\mbox{Aut}}

\def\cl{\mbox{Cl}}

\def\diag{\mbox{diag}}

\def\ker{\mbox{Ker}}

\def\max{\mbox{max}}

\def\gl{\mbox{GL}}

\def\G{{\mathbf G}}

\def\V{{\bf V}}
\def\W{{\bf W}}

\def\Z{\mathbb{Z}}

\def\p{\varphi}

\def\inv{^{-1}}
\def\la{\langle}
\def\ra{\rangle}

\def\bi{\begin{itemize}}
\def\ei{\end{itemize}}
\def\beq{\begin{equation}}
\def\eeq{\end{equation}}

\def\ds{\displaystyle}
\def\xr{\xrightarrow}

\newtheorem{T}{Theorem}[section]
\newcommand{\bt}{\begin{T}}
\newcommand{\et}{\end{T}}
\newcommand{\ftd}{$\square$\end{T}}

\newtheorem{Proposition}[T]{Proposition}
\newcommand{\bp}{\begin{Proposition}}
\newcommand{\ep}{\end{Proposition}}
\newcommand{\fpd}{$\square$\end{Proposition}}

\newtheorem{Lemma}[T]{Lemma}
\newcommand{\bl}{\begin{Lemma}}
\newcommand{\el}{\end{Lemma}}
\newcommand{\fld}{$\square$\end{Lemma}}

\newtheorem{Corol}[T]{Corollary}
\newcommand{\bc}{\begin{Corol}}
\newcommand{\ec}{\end{Corol}}
\newcommand{\fcd}{$\square$\end{Corol}}

\newtheorem{Result}[T]{Result}
\newcommand{\br}{\begin{Result}}
\newcommand{\er}{\end{Result}}
\newcommand{\frd}{$\square$\end{Result}}

\newtheorem{Example}[T]{Example}
\newcommand{\be}{\begin{Example}}
\newcommand{\ee}{\end{Example}}

\newtheorem{Problem}[T]{Problem}
\newcommand{\bq}{\begin{Problem}}
\newcommand{\eq}{\end{Problem}}

\newtheorem{Remark}[T]{Remark}
\newcommand{\brem}{\begin{Remark}}
\newcommand{\erem}{\end{Remark}}

\newtheorem{Conjecture}[T]{Conjecture}
\newcommand{\bcon}{\begin{Conjecture}}
\newcommand{\econ}{\end{Conjecture}}

\newcommand{\proof}
   {\par\medbreak\noindent{\bf Proof}.\enspace}

\newcommand{\qed}{$\Box$
\par\bigbreak}

\title{On the pseudovariety of groups $\mathbf{U} = \ds\bigvee_{p \in \mathbb{P}} \ab(p) \ast \ab(p-1)$}
\author{{\bf Claude Marion, Pedro V. Silva, Gareth Tracey}}

\date{\today}

\begin{document}

\maketitle

\begin{center}\small
2020 Mathematics Subject Classification: 20E05, 20E10, 20F10, 20F16

\bigskip

Keywords: subgroups of the free group, profinite topology, semidirect product
\end{center}

\abstract{We introduce the pseudovariety of finite groups $\mathbf{U} = \ds\bigvee_{p \in \mathbb{P}} \ab(p) \ast \ab(p-1)$, where $\mathbb{P}$ is the set of all primes. We show that $\mathbf{U}$ consists of all finite supersolvable groups with elementary abelian derived subgroup and abelian Sylow subgroups, being therefore decidable. We prove that it is decidable whether or not a finitely generated subgroup of a free group is closed or dense for the pro-$\U$ topology.
We consider also the pseudovariety of finite groups $\ab(p) \ast \ab(d)$ (where $p$ is a prime and $d$ divides $p-1$). We study the pro-$(\ab(p) \ast \ab(d))$ topology on a free group and construct the unique generator of minimum size of the pseudovariety $\ab(p) \ast \ab(d)$.
Finally, we prove that the variety of groups generated by $\U$ is the variety of all metabelian groups, obtaining also results on the varieties generated by a Baumslag-Solitar group of the form $BS(1,q)$ for $q$ prime.}

\section{Introduction}

The classification of finite groups is usually interpreted as the classification of finite simple groups, but universal algebra provides an alternative approach through the concept of a pseudovariety (a class of finite algebras closed under taking subalgebras, homomorphic images and finitary direct products). This is common practice with more general classes of algebras such as semigroups or monoids, where the role played by simple groups has no equivalent.

One of the interesting features of pseudovarieties of finite groups is that they induce a \linebreak
(pseudometrizable) topology on any group $G$. Given such a pseudovariety $\V$, the pro-$\V$ topology on $G$ is the initial topology with respect to all homomorphisms 
$G \to H \in \V$, where $H$ is endowed with the discrete topology. If ${\bf G}$ denotes the pseudovariety of all finite groups, the pro-${\bf G}$ topology is known as the profinite topology.

The profinite topology was introduced by Marshall Hall in \cite{Hal} and he proved in \cite[Theorem 5.1]{Hal2} that every finitely generated subgroup of a free group is closed for the profinite topology.

Over the years, other pseudovarieties $\V$ were considered, and the following decidability questions became objects of study for an arbitrary finitely generated subgroup $H$ of a free group $F$ (assuming that $F$ is a free group over $A$ and we are given a set of generators of $H$ expressed as reduced words over $A\cup A^{-1}$):
\bi
\item
Can we decide whether $H$ is closed for the pro-$\V$ topology?
\item
Can we decide whether $H$ is dense for the pro-$\V$ topology?
\item
Does the pro-$\mathbf{V}$ closure ${\rm Cl}_{\mathbf{V}}(H)$ of $H$ have decidable membership problem?
\item
Can we decide whether $\mathrm{Cl}_{\V}(H)$ is finitely generated and can we compute a basis in the affirmative case?
\ei

In \cite{RZ}, Ribes and Zalesski\v\i \, answered all these questions positively for the pseudovariety ${\bf G}_p$ of all finite $p$-groups, for an arbitrary prime $p$. 
In \cite{MSW}, Margolis, Sapir and Weil dealt successfully with the pseudovariety {\bf N} of all finite nilpotent groups (the case where $F$ is of infinite rank follows from \cite[Corollary 2.4]{MSTc}). The cases of the pseudovariety {\bf Ab} of all finite abelian groups and the pseudovariety {\bf M} of all finite metabelian groups are settled in \cite{MSTm}. It is understood that the case of extension-closed pseudovarieties is in general more favourable. However, the most famous open problems concern the extension-closed pseudovariety {\bf S} of finite solvable groups. 

Among the most important pseudovarieties lying between {\bf N} and {\bf S}, we find the pseudovariety {\bf Su} of finite supersolvable groups. 
A finite group is {\em supersolvable} if its chief factors are all cyclic. Equivalently, a finite group is supersolvable if every of its maximal subgroups is of prime index. 
In \cite{AS}, Auinger and Steinberg studied {\bf Su} and its Hall subpseudovarieties, and this motivated us to study the pro-{\bf Su} topology in \cite{MSTs}. In particular, we related the property of being dense for the pro-{\bf Su} topology with groups of the form $C_p \rtimes C_d$, where $p$ is a prime and $d$ divides $p-1$. 
In \cite{MSTc} we prove that cyclic subgroups of a free group of finite rank are $\mathbf{N}$-closed and consequently $\mathbf{Su}$-closed and $\mathbf{S}$-closed.
On the other hand, Auinger and Steinberg proved in \cite[Corollary 2.7]{AS} the decomposition
$${\bf Su} = \bigvee_{p \in \mathbb{P}}{\bf G}_p \ast \ab(p-1),$$
where $\mathbb{P}$ denotes the set of all primes, $\ab(n)$ denotes the class of finite groups whose exponent divides $n$ and $\ast$ is the semidirect product of pseudovarieties.

It seemed thus a natural idea to study the pseudovariety
$$\U = \bigvee_{p \in \mathbb{P}} \ab(p) \ast \ab(p-1).$$
We note that the aforementioned semidirect products $C_p \rtimes C_d$, which play an important role in the pro-{\bf Su} topology, belong to the pseudovariety $\ab(p) \ast \ab(p-1)$.

We give now a brief description of the structure of the paper. After introducing basic concepts and fixing notation in Section 2, we recall in Section 3 some known facts on the pseudovarieties $\ab(p) \ast \ab(d)$ (where $p \in \mathbb{P}$ and $d$ divides $p-1$) and compute the closure of $H \leq_{f.g.} F_n$ for the pro-$(\ab(p) \ast \ab(d))$ topology.
In Section 4, we compute the unique generator of minimum size (up to isomorphism) for each pseudovariety $\ab(p) \ast \ab(d)$.

In Section 5, we start the study of $\U$. The main result states that a finite group $G$ belongs to $\U$ if and only if it is supersolvable, its derived subgroup is elementary abelian and has abelian Sylow subgroups. As a corollary, $\U$ has decidable membership problem (i.e., there is an algorithm which decides whether or not a given arbitrary finite group belongs to $\U$).

In Section 6, we prove that the variety of groups generated by $\U$ is the variety of all metabelian groups, obtaining also results on the varieties generated by a Baumslag-Solitar group of the form $BS(1,q)$ for $q$ prime. The stronger version of some results depends on a generalized form of the Riemann hypothesis.
An intriguing conjecture is also suggested, involving the variety generated by families of groups of the form $C_p \rtimes C_{p-1}$ ($p$ prime).

We note that $\U$ is not extension-closed: it is easy to check that $A_4 \notin \U$, but $[A_4,A_4] \cong C_2^2 \in \U$ and $A_4/[A_4,A_4] \cong C_3 \in \U$. In Section 7, we show that $\U \subset {\bf Su} \cap \M$
and study the pro-$\U$ topology on $F_n$. In particular, we provide a noneffective description of the pro-$\U$ closure of $H \leq_{f.g.} F_n$ and show that is decidable whether or not $H$ is closed (or dense) for the pro-$\U$ topology on $F_n$. If $H \leq_{f.g.} F$ and $F$ is a free group of infinite rank, we show that $H$ is neither closed nor dense.

Finally, we produce a finitely generated subgroup of the free object of the variety generated by $\U$ which is not $\U$-closed, unlike the abelian and metabelian cases.

\section{Preliminaries}
\label{preli}

Given a group $G$ and elements $g, h \in G$, we use the notation $g^h=hgh^{-1}$ for the conjugate of $g$ by $h$ and $[g,h] = ghg\inv h\inv$ for the commutator of $g$ and $h$, and denote the derived subgroup of $G$ by $[G,G]$. The second derived subgroup $[[G,G],[G,G]]$ is denoted by $G^{(2)}$. 
Given subgroups $H$ and $K$ of $G$, we let $[H,K]=\langle [h,k]: h\in H, k \in K\rangle$. If $H$ and $K$ are normal in $G$ then so is $[H,K]$.
Given $H \unlhd G$ and $k \geq 1$, we denote by ${\rm pow}_k(H)$ the (normal) subgroup of $G$ generated by all powers of the form $u^k$ with $u \in H$. On the other hand, given a group $H$, $H^k$ denotes the direct product of $k$ copies of $H$.
Given a subset $X$ of a group $G$, we denote by $\la\la X\ra\ra_G$ the normal 
closure of $X$ in $G$, i.e. the smallest normal subgroup of $G$ containing $X$. We denote by $o(g)$ the order of an element $g$ of a group. 

Given a positive integer $m$, we let $C_m$ or $\mathbb{Z}/m\mathbb{Z}$   denote the cyclic group of order $m$. When the context is clear, we sometimes view  $\mathbb{Z}/m\mathbb{Z}$ as a ring (namely the ring of integers modulo $m$),  and denote by  $(\mathbb{Z}/m\mathbb{Z})^*$ the group of units of $\mathbb{Z}/m\mathbb{Z}$. When $p$ is a prime, we also write  $\mathbb{F}_p$ for $\mathbb{Z}/p\mathbb{Z}$.

Let $p$ be a prime and let $d > 1$ be a divisor of $p-1$. It is well known that $(\mathbb{Z}/p\mathbb{Z})^*$ is a cyclic group. Hence $(\mathbb{Z}/p\mathbb{Z})^*$ has a unique subgroup of order $d$, which is itself cyclic. 
Given a positive integer $k$, we let $[k]=\{1,\dots,k\}\subset \mathbb{N}$.
We shall use the notation
$$Q_{p,d} = \{ q \in [p-1] \mid q^d \equiv 1\,({\rm mod}\, p)\} \mbox{ and } Q'_{p,d} = \{ q \in [p-1] \mid q\mbox{ has order $d$ in }\mathbb{F}_p^*\}.$$

Given a ring $R$ and elements $k_1,\dots,k_n\in R$, we let ${\rm diag}(k_1,\dots,k_n)$ be the matrix having entries $k_1,\dots,k_n$ in the main diagonal. Given a direct product $R = R_1 \times \ldots \times R_k$ of rings with unity and $i \in [k]$, we use the notation $e_i$ for the element of $R$ having 1 in the $i$th component and 0 everywhere else.
Let $1^k$ be shorthand for a sequence $1,1,\ldots,1$ of $k$ 1's. 

If $A$ is a set, we set $F_A$ to be the free group on $A$. The length of an element $u \in F_A$ is the length of the reduced form of $u$ and is denoted by $|u|$. 
Also if $H$ is a finitely generated subgroup of $F_A$, we write $H\leq_{f.g.} F_A$. We denote by $F_n$ the free group of rank $n \in \N$ and denote by $A$ or $A_n$ a fixed basis of $F_n$.

In this paper a (finite) {\em automaton} is a structure of the form $\A = (A, Q, q_0, E)$, where the {\em alphabet} $A$ and the {\em vertex set} $Q$ are (finite) sets, $q_0 \in Q$ is the {\em basepoint} and $E \subseteq Q \times A \times Q$ is the {\em edge set}.

Given a subgroup $H$ of a free group $F_A$, we define the {\em Schreier automaton} $S(H) = (A \cup A\inv,Q,H,E)$ by
$$Q = \{ Hu \mid u \in F_A\},\quad E = \{ (Hu,a,Hua) \mid u\in F_A, a \in A\cup A\inv\}.$$
If $H$ is finitely generated, the subautomaton $\A(H)$ induced by the set of vertices of the form $Hu$, where $u$ is a prefix of some $h \in H$, is a finite automaton known as the {\em Stallings automaton} of $H$. In \cite{Sta}, Stallings introduced an algorithm to compute $\A(H)$ from any finite generating set of $H$ through the {\em folding} operation), and provided several important applications of the concept. For the basic properties of Stallings automata, the reader is referred to \cite[Section 3]{BS}.

A {\em pseudovariety of finite groups} is a class of finite groups closed under taking subgroups, homomorphic images and finitary direct products. Since no other pseudovarieties occur in the paper, we will just write {\em pseudovariety}. For the general theory of pseudovarieties, the reader is referred to \cite{RS}.

Classical examples include the pseudovariety $\ab$ of all finite abelian groups and the Burnside pseudovarieties $\bu(n)$. For every $n \geq 1$, $\bu(n)$ denotes the class of all finite groups $G$ such that $g^n = 1$ for every $g \in G$. Since the intersection of pseudovarieties is still a pseudovariety, then
$$\ab(n) = \ab \cap \bu(n)$$
is also a pseudovariety for every $n \geq 1$. In other words $\ab(n)$ is the pseudovariety of all finite abelian groups of exponent dividing $n$.

If $\V$ and $\W$ are pseudovarieties, the pseudovariety $\V$-by-$\W$ contains all finite groups $G$ admitting a normal subgroup $N \in \V$ such that $G/N \in \W$. As noted by Auinger and Steinberg in \cite{AS}, it follows from the Kalu\v{z}nin-Krasner theorem \cite[Theorem 22.21]{Neu} that $\mathbf{V}$-by-$\mathbf{W}=\mathbf{V}*\mathbf{W}$, i.e. the pseudovariety of finite groups containing all subgroups of semidirect products $V\rtimes W$ where $V\in \mathbf{V}$ and $W\in \mathbf{W}$.

A {\em variety of groups} is a class of groups closed under taking subgroups, homomorphic images and arbitrary direct products. By Birkhoff's Theorem, a class of groups constitutes a variety if and only if it is the class of all groups satisfying a certain set of group identities. If $X$ denotes a countable alphabet, a {\em group identity} is a formal equality of the type $u = v$, where $u,v \in F_X$. A group $G$ satisfies this identity if $\p(u) = \p(v)$ for every homomorphism $\p:F_X \to G$. 
We then write $G \models u = v$. Given a set of identities $\Sigma \subseteq F_X \times F_X$, we denote by $[\Sigma]$ the variety of all groups satisfying each one of the identities in $\Sigma$.

We consider two important examples of varieties in this paper. One is the variety of all abelian groups, denoted by $\cal{AB}$ and defined by the single identity $xy = yx$. The other is the variety of all metabelian groups, denoted by $\cal{M}$. A group $G$ is {\em metabelian} if $[G,G]$ is abelian, hence $\cal{M}$ is defined by the single identity $[x_1,x_2][x_3,x_4] = [x_3,x_4][x_1,x_2]$.

Given a class $C$ of groups, we denote by ${\cal{V}}(C)$ the variety of groups {\em generated} by $C$. This is the smallest variety containing $C$ and consists of all homomorphic images of subgroups of direct products of groups in $C$. 

Let ${\cal{V}} = [\Sigma]$.
Given an alphabet $A$, let
$$F_A^{{\cal{V}}} = \la\la \p(uv\inv) \mid (u,v) \in \Sigma,\, \p:F_X \to F_A\mbox{ homomorphism}\ra\ra_{F_A} \unlhd F_A.$$
It follows easily from the definitions and the universal property of $F_A$ that $F_A({\cal{V}}) = F_A/F_A^{{\cal{V}}}$ is the {\em free object} of $\cal{V}$ on $A$. 
For details on varieties of groups, the reader is referred to \cite{Neu}.

Now ${\cal{V}}^f\, =\, {\cal{V}} \cap \G$ is clearly a pseudovariety, known as the {\em finite trace of} $\cal{V}$. In view of Birkhoff's Theorem, a pseudovariety $\V$ is of the form $\mathbf{V}=\mathcal{V}^f$ if and only if it is {\em equational}, that is, there exists a set $\Sigma$ of group identities such that $\V = [\Sigma] \cap \G$. Note that arbitrary pseudovarieties require in general pseudoidentities (in view of Reiterman's Theorem). The reader is referred to \cite[Section 7.2]{RS} for further details.

Let $\V$ denote a pseudovariety and let $\Sigma$ be a set of identities. We write $\V \models \Sigma$ if $G \models \sigma$ for all $G \in \V$ and $\sigma \in \Sigma$. Write $[[\Sigma ]] = [\Sigma ]^f$.

Given a finite group $G$, we denote by $\langle G\rangle$ the pseudovariety generated by $G$ (consisting of all homomorphic images of subgroups of direct powers of the form $G^n$). Such a pseudovariety is said to be {\em finitely generated}. Note that the smallest pseudovariety containing the finite groups $G_1, \ldots, G_n$ is the pseudovariety generated by $G_1 \times \ldots \times G_n$ (this explains the terminology).

Every finitely generated pseudovariety $\V$ admits {\em free objects}. That is, given any finite set $A$, there exists some  $F_A(\V) \in \V$ and some mapping $\iota:A \to F_A(\V)$ such that, for every mapping $\p:A \to G \in \V$, there exists some (unique) homomorphism $\Phi:F_A(\V) \to G$ such that $\Phi \circ \iota = \p$. Moreover, $F_A(\V)$ is unique up to isomorphism. On the other hand, it follows from 
\cite{Rei} that $\V = [[\Sigma ]]$ for some set of identities $\Sigma$.

Given a pseudovariety $\mathbf{V}$, where we consider finite groups endowed with the discrete topology, the pro-$\mathbf{V}$ topology on a group $G$ is defined as the coarsest topology which makes all morphisms from $G$ into elements of $\mathbf{V}$ continuous. Equivalently, $G$ is a topological group where the normal subgroups $K$ of $G$ such that $G/K\in \mathbf{V}$ form a basis of neighbourhoods of the identity. 
By \cite[Proposition 1.2]{MSW}, a subgroup $H$ of $G$ is pro-$\mathbf{V}$ open if and only if it is pro-$\mathbf{V}$ clopen, and if and only if $G/\textrm{Core}_G(H)$ belongs to $\mathbf{V}$ (the core $\textrm{Core}_G(H)$ of $H$ in $G$ is the largest normal subgroup of $G$ contained in $H$ and is equal to $\bigcap_{g\in G} g^{-1}Hg$). Note that $H$ is $\mathbf{V}$-closed if and only, for every $g\in G\setminus H$, there exists some $\mathbf{V}$-clopen $K\leq G$ such that $H \leq K$ and $g\not \in K$. Moreover, a subgroup $H$ of $G$ is pro-$\mathbf{V}$ dense if and only if $HN=G$ for every normal subgroup $N$ of $G$ such that $G/N\in \mathbf{V}$.

For a topological property $\mathcal{P}$ and a subset $S$ of $G$, we say that $S$ is $\mathbf{V}$-$\mathcal{P}$ if $S$ has property $\mathcal{P}$ in the pro-$\mathbf{V}$ topology on $G$. 
Given $S \subseteq G$, we also denote by  ${\rm{Cl}}_{\mathbf{V}}^G(S)$ the $\mathbf{V}$-closure of $S$ in $G$ and we omit the superscript when no confusion is possible.
In this paper we mainly consider the case $G = F_n$. If $H \leq G$, then also ${\rm{Cl}}_{\mathbf{V}}^G(H) \leq G$ \cite[Proposition 1.3]{MSW}.

Suppose that $\V$ and $\mathbf{W}$ are pseudovarieties of finite groups. Then:
\beq
\label{icl}
\mbox{If $\mathbf{W} \subseteq \V$, then ${\rm{Cl}}_{\mathbf{V}}^G(S) \subseteq {\rm{Cl}}_{\mathbf{W}}^G(S)$ for every $S \subseteq G$.}
\eeq
This follows from the pro-$\mathbf{W}$ topology being coarser than the pro-$\mathbf{V}$ topology on $G$.

Throughout the paper, unless otherwise stated, we fix a prime $p > 2$ and a divisor $d > 1$ of $p-1$.

\section{The pseudovarieties $\ab(p) \ast \ab(d)$}

Let $p>2$ be a prime and let $d>1$ be a divisor of $p-1$. 

By a theorem of \v{S}mel'kin \cite{Sme} (see also \cite[Theorem 24.64]{Neu}), the pseudovariety $\ab(p) \ast \ab(d)$ is finitely generated and admits therefore free objects. 
For every $n \geq 1$, let $F_n(p,d)$ denote the free object of $\ab(p) \ast \ab(d)$
on a set of size $n$. By a theorem of Higman \cite{Hig} (see also \cite[Item 24.65]{Neu}), we have
$$\ab(p) \ast \ab(d) = \la F_2(p,d)\ra.$$
How big is $F_2(p,d)$? 
For all $m,n \geq 1$, we have $\ab(m) = \la C_m\ra$ and it follows easily that 
$$F_n(\ab(m)) \cong C_m^n.$$
Now it follows from \cite[Corollary 21.13]{Neu} that $F_n(p,d)$ admits a semidirect product decomposition
$$ F_n(p,d) = C_p^{(n-1)d^n+1} \rtimes C_d^n.$$
In particular, $|F_2(p,d)| = p(d^2+1)d^2$. 

Since $\ab(m) = \bigl[\bigl[ [x,y],x^m \bigr]\bigr]$ for every $m \geq 1$, we have
$$F_n(\ab(m)) \cong F_n/\la\la [F_n,F_n] \cup {\rm pow}_m(F_n)\ra\ra_{F_n}.$$
Write $K_{n,m} = [F_n,F_n]{\rm pow}_m(F_n)$. It is routine to check that $K_{n,m} \unlhd F_n$, hence $\la\la [F_n,F_n] \cup {\rm pow}_m(F_n)\ra\ra_{F_n} = K_{n,m}$ and so
$$F_n(\ab(m)) \cong F_n/K_{n,m}.$$
Write
$$L_{n,p,d} = [K_{n,d},K_{n,d}]{\rm pow}_p(K_{n,d}).$$
It follows from \cite[Item 21.12]{Neu} that 
\beq
\label{iden}
\ab(p) \ast \ab(d) = \Bigl[\Bigl[ w = 1 \; (w \in \bigcup_{n \geq 1} L_{n,p,d}) \Bigr]\Bigr].
\eeq
Now we can derive easily the following alternative characterization:

\bp
\label{altfree}
Let $p > 2$ be a prime and let $d>1$ be a divisor of $p-1$. If $n \geq 1$ then $F_n(p,d) \cong F_n/L_{n,p,d}$.
\ep

\proof
In view of (\ref{iden}), it suffices to show that $F_n/L_{n,p,d} \in \ab(p) \ast \ab(d)$, which amounts to proving that 
$F_n/L_{n,p,d} \models w = 1$ for all $w \in L_{m,p,d}$ and $m \geq 1$. It is therefore enough to show that
$\p(w) \in L_{n,p,d}$ for every homomorphism $\p:F_m \to F_n$ and every $w \in L_{m,p,d}$. This can be further reduced to proving that $\p(u) \in K_{n,d}$ for every homomorphism $\p:F_m \to F_n$ and every $u \in K_{m,d}$, which is straightforward.
\qed

We can now compute the $(\ab(p) \ast \ab(d))$-closure of a subgroup of $F_n$. 

\bp
\label{clop}
Let $p>2$ be a prime and let $d>1$ be a divisor of $p-1$. Let $H \leq F_n$ and let ${\rm Cl}(H)$ denote the $(\ab(p) \ast \ab(d))$-closure of $H$ in $F_n$. Then:
\bi
\item[(i)]
${\rm Cl}(H) = HL_{n,p,d}$;
\item[(ii)]
${\rm Cl}(H)$ has finite index in $F_n$;
\item[(iii)]
${\rm Cl}(H)$ is finitely generated.
\ei
\ep

\proof
(i) We show that
\beq
\label{clop1}
\cl(1) = L_{n,p,d}.
\eeq

Since $L_{n,p,d} \unlhd F_n$ and $F_n/L_{n,p,d} \cong F_n(p,d) \in \ab(p) \ast \ab(d)$ by Proposition \ref{altfree}, it follows that $L_{n,p,d}$ is a clopen subgroup of $F_n$. Thus $\cl(1) \leq L_{n,p,d}$. On the other hand, $\cl(1)$ is the intersection of all the open subgroups of $F_n$ by \cite[Proposition 1.3]{MSW}, hence it suffices to show that $L_{n,p,d} \leq K$ for every open $K \leq F_n$.

Given such $K$, we have 
$F_n/\textrm{Core}_{F_n}(K) \in \ab(p) \ast \ab(d)$. Since we may assume that $F_n/L_{n,p,d} = F_n(p,d)$, the canonical mapping of the basis alphabet $A_n$ to  
$F_n/\textrm{Core}_{F_n}(K)$ induces a homomorphism from $F_n/L_{n,p,d}$ onto $F_n/\textrm{Core}_{F_n}(K)$, yielding $L_{n,p,d} \leq \textrm{Core}_{F_n}(K) \leq K$. Therefore (\ref{clop1}) holds.

Since $H,  L_{n,p,d} \leq \cl(H)$, we get
$HL_{n,p,d} \leq \cl(H)$. Thus it suffices to show that $HL_{n,p,d}$ is closed. Now $L_{n,p,d} \unlhd F_n$ yields
$L_{n,p,d} \leq \textrm{Core}_{F_n}(HL_{n,p,d})$. Hence $F_n/\textrm{Core}_{F_n}(HL_{n,p,d})$ is a quotient of $F_n/L_{n,p,d} \in \ab(p) \ast \ab(d)$ and so $F_n/{\rm Core}_{F_n}(HL_{n,p,d}) \in \mathbf{Ab}(p)*\mathbf{Ab}(d)$. Hence $HL_{n,p,d}$ is open and so closed, as required.

(ii) In view of Proposition \ref{altfree}, we have $[F_n:L_{n,p,d}] < \infty$. Now it follows from part (i) that $[F_n:\cl(H)] < \infty$.

(iii) By part (ii), since $F_n$ is finitely generated.
\qed

\bc
\label{clopc}
Let $p>2$ be a prime and let $d>1$ be a divisor of $p-1$. Let $H \leq_{f.g.} F_n$. Then the $(\ab(p) \ast \ab(d))$-closure of $H$ in $F_n$ is computable from a finite generating subset of $H$.
\ec

\proof
Let ${\rm Cl}(H)$ denote the $(\ab(p) \ast \ab(d))$-closure of $H$ in $F_n$. Recall that $L_{n,p,d} \unlhd F_n$. Since the structure of the finite group $F_n(p,d) \cong F_n/L_{n,p,d}$ is perfectly known, we can construct the Cayley graph of $F_n/L_{n,p,d}$ with respect to $A_n$, which is also the Stallings automaton $\A(L_{n,p,d})$ of $L_{n,p,d}$ with respect to the basis $A_n$  of $F_n$. And we can also construct $\A(H)$ since we are given a finite generating set of $H$.
By Proposition \ref{clop}, we have
$$\cl(H) = HL_{n,p,d}  = H \vee L_{n,p,d}.$$
Now it follows from Stallings' algorithm \cite[Section 3]{BS} that $\A(H \vee L_{n,p,d})$ is obtained by identifying the basepoints of $\A(L_{n,p,d})$ and $\A(H)$ and performing a complete folding. Therefore $\A(\cl(H))$ is effectively constructible, which is equivalent to saying that $\cl(H)$ is computable.
\qed

\section{A generator of minimum size for $\ab(p) \ast \ab(d)$}

Let $p>2$ be a prime and let $d>1$ be a divisor of $p-1$. 
Since $p$ and $d$ are coprime, it follows from Schur-Zassenhaus theorem \cite[Theorem 7.41]{Rot} that the pseudovariety $\ab(p)$-by-$\ab(d) = \ab(p) \ast \ab(d)$ consists of all semidirect products of the form $N \rtimes H$ with $N \in \ab(p)$ and $H \in \ab(d)$.

The (outer) semidirect product of two groups $N$ and $H$ is defined by a group action 
$$\begin{array}{rcl}
\p:H&\to&\aut(N)\\
h&\mapsto&\p_h
\end{array}$$ 
which induces a group operation on $N \times H$ by
$$(u,h)(u',h') = (u\cdot \p_h(u'),hh').$$
This group is denoted by $N \rtimes_{\p} H$. Note that:
\bi
\item
$N \cong N \times \{ 1\} \unlhd N \rtimes_{\p} H$;
\item
$H \cong \{ 1\} \times H \leq N \rtimes_{\p} H$;
\item
$(1,h)(u,1)(1,h)\inv = (\p_h(u),h)(1,h\inv) = (\p_h(u),1)$ for all $u \in N$ and $h \in H$.
\ei
It is well known how to build a presentation for $N \rtimes_{\p} H$ from presentations for $N$ and $H$. Assume that $\la X \mid R\ra$ is a presentation for $N$ and $\la Y \mid S\ra$ is a presentation for $H$, with $X \cap Y = \emptyset$. For all $x \in X$ and $y \in Y$, fix a word $\oo{\p}_y(x) \in (X \cup X\inv)^*$ (the free monoid on $X \cup X\inv$) representing $\p_y(x) \in N$. Then 
\beq
\label{psp}
\la X \cup Y \mid R \cup S, \, yxy\inv = \oo{\p}_y(x)\; (x \in X,\, y \in Y) \ra.
\eeq
is a presentation of $N \rtimes_{\p} H$.

We claim that we can restrict the type of semidirect product required to produce all the groups in $\ab(p) \ast \ab(d)$.
Indeed, assume that $N \in \ab(p)$. Then $N$ is a direct product of $n$ copies of $\Z/p\Z$ for some $n \geq 1$ and so $\aut(N)$ is isomorphic to the group of linear automorphisms of a vector space of rank $n$ over the field ${\mathbb{F}}_p$, which in turn is isomorphic to $\gl_n({\mathbb{F}}_p)$. 

Since we are writing functions on the left, we must be careful on defining this isomorphism. Let $\{ e_1,\ldots,e_n\}$ denote the canonical basis of $\mathbb{F}_p^n$. Given $\theta \in \aut(N)$, let $M_{\theta}\in \gl_n({\mathbb{F}}_p)$ be the matrix where the $i$th row vector is given by the coordinates of $\theta(e_i)$ in the canonical basis. Then 
$$\aut(N) \to \gl_n({\mathbb{F}}_p): \theta \mapsto M_{\theta}$$
is really an anti-isomorphism. To get a true isomorphism, we must define instead
$$\begin{array}{rcl}
\Lambda:\aut(N)&\to&\gl_n({\mathbb{F}}_p)\\
\theta&\mapsto&M_{\theta\inv}
\end{array}$$

We say that $\p:H \to \aut(N)$ is a {\em diagonalizable} action if there exists some $P \in \gl_n({\mathbb{F}}_p)$ such that $P\inv M_{\p_h}P$ is a diagonal matrix for every $h \in H$. 

\bl
\label{diagfp}
Let $n$ be a positive integer and let $p$ be a prime. If $M \in {\rm GL}_n(\mathbb{F}_p)$ and $M^{p-1}=1$ then $M$ is diagonalizable over $\mathbb{F}_p$.
\el

\proof
 Let $K$ be an algebraic closure of $\mathbb{F}_p$. Note  that $\mathbb{F}_p$ is equal to $\{a \in K: a^p=a\}$. Indeed $\{a \in K: a^p=a\}$ is a subfield of $K$ of size $p$ containing $\mathbb{F}_p$ (and so must be $\mathbb{F}_p$). Consider $M$ as an element of ${\rm GL}_n(K)$. By Jordan decomposition, there exist unique elements $s$ and $u$ in ${\rm GL}_n(K)$ such that $M=su=us$ where $s$ is semisimple (i.e. diagonalizable) and $u$ is unipotent (i.e. all the eigenvalues of $u$ are equal to 1). Note that since ${\rm char}(K)>0$, an element $u$ of ${\rm GL}_n(K)$ is unipotent if and only if $u$ has order a power of $p$. 

Now since $M$ has order dividing $p-1$, then $u^{p-1} = s^{1-p}$ is both unipotent and diagonalizable, hence $u^{p-1} = 1$ and so $u^p = u$. Since $u$ has order a power of $p$, we get that $u=1$ and $M=s$, that is $M$ is diagonalizable over $K$. Hence there exist $a_1,\dots,a_n$ in $K$ and $T$ in ${\rm GL}_n(K)$ such that $T^{-1}MT={\rm diag}(a_1,\dots,a_n)$.
On one hand, we have
$$(T^{-1}MT)^p={\rm diag}(a_1^p,\dots,a_n^p).$$

On the other hand, we have
\begin{eqnarray*}
(T^{-1}MT)^p&= & (T^{-1}MT)^{p-1} (T^{-1}MT)\\
                     &=&T^{-1}M^{p-1}T (T^{-1}MT)\\
                     &=& T^{-1}MT\\
                     &=&{\rm diag}(a_1,\dots,a_n).
\end{eqnarray*}

We therefore obtain $a_i^p=a_i$ and so $a_i\in \mathbb{F}_p$ for $1\leq i \leq n$. In particular ${\rm diag}(a_1,\dots,a_n)$ is in ${\rm GL}_n(\mathbb{F}_p)$. Now $M$ and ${\rm diag}(a_1,\dots,a_n)$ are both in ${\rm GL}_n(\mathbb{F}_p)$ and are conjugate in ${\rm GL}_n(K)$.  It follows that  $M$ and ${\rm diag}(a_1,\dots,a_n)$ are conjugate in ${\rm GL}_n(\mathbb{F}_p)$ (see \cite[Theorem 4]{Con} for a proof, which can be derived from the uniqueness of the Frobenius normal form over the smaller field). Hence $M$ is diagonalizable over $\mathbb{F}_p$.
\qed

\bl
\label{alwdia}
Let $p>2$ be a prime and let $d>1$ be a divisor of $p-1$. Let $N \in \ab(p)$ and $H \in \ab(d)$. Every action $\p:H \to \aut(N)$ is diagonalizable.
\el

\proof
For every $h \in H$, let $M'_h = \Lambda(\p_h) \in \gl_n({\mathbb{F}}_p)$. We have
$$(M'_h)^{p-1} = (\Lambda(\p_h))^{p-1} = (M_{\p_h\inv})^{p-1} = (M_{\p_{h\inv}})^{p-1} = M_{\p_{(h\inv)^{p-1}}} = M_1 = 1,$$
hence $M'_h$ is diagonalizable over ${\mathbb{F}}_p$ by Lemma \ref{diagfp}. On the other hand, for all $h,h' \in H$, since $H$ is abelian, we have
$$M'_hM'_{h'} = \Lambda(\p_h) \Lambda(\p_{h'}) = \Lambda(\p_h\p_{h'}) = \Lambda(\p_{hh'}) = \Lambda(\p_{h'h}) = \Lambda(\p_{h'}\p_{h}) = \Lambda(\p_{h'}) \Lambda(\p_{h}) = M'_{h'}M'_{h}.$$
Since the matrices $M'_h$ ($h \in H$) commute and are diagonalizable over $\mathbb{F}_p$, they are simultaneously diagonalizable over $\mathbb{F}_p$ (see \cite[Theorem 5]{Con1} for a proof).
Thus there exists some $P \in \gl_n({\mathbb{F}}_p)$ such that $P\inv M'_hP = P\inv M_{\p_{h\inv}}P$ is a diagonal matrix for every $h \in H$. Therefore $\p$ is diagonalizable.
\qed

\bl
\label{diapr}
Let $p>2$ be a prime and let $d>1$ be a divisor of $p-1$. The following conditions are equivalent for a group $G$:
\bi
\item[(i)] $G \in \ab(p) \ast \ab(d)$;
\item[(ii)] $G$ admits a presentation of the form
\beq
\label{diapr1}
\begin{array}{l}
\la x_1,\ldots,x_n,y_1,\ldots,y_m \mid x_i^p = 1,\, y_j^{d_j} = 1, \, x_ix_{i'} = x_jx_{i'}, \, y_jy_{j'} = y_{j'}y_j,\\ \\
\hspace{4cm} y_jx_iy_j\inv = x_i^{q_{ij}}\, (i,i' \in [n],\, j,j' \in [m])\ra
\end{array}
\eeq
for some $n,m \geq 0$, divisors $d_j > 1$ of $d$ and $q_{ij} \in Q_{p,d_j}$ $(i \in [n],\, j \in [m])$.
\ei
\el

\proof
(i) $\Rw$ (ii). Assume that $G = N \rtimes_{\p} H$ with $N \in \ab(p)$ and $H \in \ab(d)$. By Lemma \ref{alwdia}, the action $\p:H \to \aut(N)$ is diagonalizable. We may assume that both $N$ and $H$ are nontrivial. Then $N \cong (\Z/p\Z)^n$ for some $n \geq 1$ and
$H \cong (\Z/d_1\Z) \times \ldots \times (\Z/d_m\Z)$ for some divisors $d_1,\ldots,d_m > 1$ of $d$. 
Let $a_i \in N$ correspond to $e_i \in (\Z/p\Z)^n$ for each $i \in [n]$, and let $y_j \in H$ correspond to $e_j \in (\Z/d_1\Z) \times \ldots \times (\Z/d_m\Z)$ for each $j\in [m]$. Then $N = \la a_1,\ldots,a_n\ra$ and $H = \la y_1,\ldots,y_m\ra$.

Since $\p$ is diagonalizable, there exists some $P \in \gl_n({\mathbb{F}}_p)$ such that $P\inv M_{\p_h}P$ is a diagonal matrix for every $h \in H$. Note that $M_{\p_h}$ is invertible since $M_{\p_h}M_{\p_{h\inv}}$ is the identity matrix.
Then there exist 
$q_{ij} \in [p-1]$ $(i \in [n],\, j \in [m])$ such that
$$P\inv M_{\p_{y_j}}P = \diag(q_{1j},\ldots,q_{nj})$$
for each $j \in [m]$.

For each $i \in [n]$, write the Kronecker vector $e_i \in \mathbb{F}_p^n$ as $(\delta_{1i},\ldots,\delta_{ni})$ and let $x_i \in N$ correspond to the $i$th column of the matrix $P = (p_{ii'})$. Then
$$\begin{array}{lll}
M_{\p_{y_j}}(p_{1i}, p_{2i}, \ldots, p_{ni})^T&=&PP\inv M_{\p_{y_j}}Pe_i^T = P\diag(q_{1j},\ldots,q_{nj})e_i^T\\ &&\\
&=&P\diag(q_{1j}\delta_{1i},\ldots,q_{nj}\delta_{ni}) = (p_{1i}q_{ij}, p_{2i}q_{ij}, \ldots, p_{ni}q_{ij})^T
\end{array}$$
and so $\p_{y_j}(x_i) = x_i^{q_{ij}}$ for all $i \in [n]$ and $j \in [m]$ (in multiplicative notation). Hence
$$(1,y_j)(x_i,1)(1,y_j)\inv = (\p_{y_j}(x_i),y_j)(1,y_j\inv) = (x_i^{q_{ij}},1).$$
We may abuse notation by identifying $x_i$ with $(x_i,1)$ and $y_j$ with $(1,y_j)$. It follows that
$$G = \la x_1,\ldots,x_n,y_1,\ldots,y_m\ra$$ and satisfies all the relations in (\ref{diapr1}). 

Given $i \in [n]$ and $j \in [m]$, it is easy to prove by induction on $k \geq 0$ that
\beq
\label{ndiapr1}
y_j^kx_iy_j^{-k} = x_i^{q_{ij}^k}
\eeq
holds for the presentation (\ref{diapr1}).
It follows that
$$x_i = y_j^{d_j}x_iy_j^{-d_j} = x_i^{q_{ij}^{d_j}},$$
hence $q_{ij}^{d_j} \equiv 1 \,({\rm mod}\, p)$ and so $q_{ij} \in Q_{p,d_j}$.

Let $G'$ be the group defined by 
(\ref{diapr1}). Since $G$ is a quotient of $G'$, it suffices to show that $|G'| \leq |G|$. Clearly, $|G| = p^nd_1\ldots d_m$. 
Now, using the relations in (\ref{diapr1}), we see that every element of $G'$ can be written in the form 
$$x_1^{r_1}\ldots x_n^{r_n}y_1^{s_1}\ldots y_m^{s_m}$$
with $r_i \in \{ 0,\ldots,p-1\}$ for $i \in [n]$ and $s_j \in \{ 0, \ldots,d_j-1\}$ for $j \in [m]$. Therefore $|G'| \leq p^nd_1\ldots d_m$ and so (\ref{diapr1}) is a presentation of $G$ as desired.

(ii) $\Rw$ (i). Straightforward. 
\qed

Given $q \in Q'_{p,d}$, we denote by $G_{p,d}^{(q)}$ the group generated by the presentation
\beq
\label{presg}
\la x,y \mid x^p = 1,\, y^d = 1, \,  yxy\inv = x^q\ra.
\eeq

\bl
\label{proot}
Let $p>2$ be a prime and let $d>1$ be a divisor of $p-1$. If $q,r\in Q'_{p,d}$, then $G_{p,d}^{(q)} \cong G_{p,d}^{(r)}$.
\el

\proof
Since ${\mathbb{F}}_p^*$ has a unique (cyclic) subgroup of order $d$, there exist $m,k \geq 1$ such that $r^m = q$ and $q^k = r$ modulo $p$. Consider $\{x,y\}$ as a basis for $F_2$. Let $\p,\psi$ be the endomorphisms of $F_2$ defined by
$$\p(x) = \psi(x) = x,\quad \p(y) = y^m,\quad \psi(y) = y^k.$$
It is routine to check that $\p$ and $\psi$ induce mutually inverse isomorphisms $\oo{\p}:G_{p,d}^{(q)} \to G_{p,d}^{(r)}$ and  $\oo{\psi}:G_{p,d}^{(r)} \to G_{p,d}^{(q)}$. For instance,
$$\oo{\p}(yxy\inv) = y^mxy^{-m} = x^{r^m} = x^q = \oo{\p}(x^q),$$
while $r^{mk} \equiv  r\,({\rm mod}\, p)$ yields $mk \equiv 1\,({\rm mod}\, d)$ and consequently
$$\oo{\psi}(\oo{\p}(y)) = y^{mk} = y.$$
Therefore $G_{p,d}^{(q)} \cong G_{p,d}^{(r)}$.
\qed 

Therefore we can omit the exponent in $G_{p,d}^{(q)}$.

\bp
\label{figen}
Let $p>2$ be a prime and let $d>1$ be a divisor of $p-1$. Then $\ab(p) \ast \ab(d) = \la G_{p,d}\ra$. Moreover, $G_{p,d}$ is up to isomorphism the unique generator of minimum size of $\ab(p) \ast \ab(d)$.
\ep

\proof
By Lemma \ref{diapr}, we have $G_{p,d} \in \ab(p) \ast \ab(d)$, hence $\la G_{p,d}\ra \subseteq 
\ab(p) \ast \ab(d)$.

Conversely, let $G \in \ab(p) \ast \ab(d)$. By Lemma \ref{diapr}, we may assume that $G$ is defined by a presentation of the form (\ref{diapr1}). For every $i \in [n]$, let $G_i$ be the group defined by the presentation
$$\la x_i,y_1,\ldots,y_m \mid x_i^p = 1,\, y_j^{d_j} = 1, \, y_jy_{j'} = y_{j'}y_j,\, y_jx_iy_j\inv = x_i^{q_{ij}}\, (j,j' \in [m])\ra.$$
We define a homomorphism
$\theta:G\to G_1 \times \ldots \times G_n$ by
$$\theta(x_i) = (1^{i-1},x_i,1^{n-i})\quad \mbox{and}\quad \theta(y_j) = (y_j,\ldots,y_j)\quad (i \in [n],\; j \in [m]).$$
It is straightforward to check that $\theta$ is a well-defined homomorphism. For instance,
$$\begin{array}{lll}
\theta(y_jx_iy_j\inv)&=&(y_j,\ldots,y_j)(1^{i-1},x_i,1^{n-i})(y_j\inv ,\ldots,y_j\inv) = (1^{i-1},y_jx_iy_j\inv,1^{n-i})\\ &&\\
&=& (1^{i-1},x_i^{q_{ij}},1^{n-i}) = \theta(x_i^{q_{ij}}).
\end{array}$$

Let $g \in \ker(\theta)$. Then we may write
$$g = x_1^{r_1}\ldots x_n^{r_n}y_1^{s_1}\ldots y_m^{s_m}$$
for some $r_i \in \{ 0,\ldots,p-1\}$ and $s_j \in \{ 0, \ldots, d_j-1\}$ $(i \in [n],\, j \in [m])$. Hence
$$(1^n) = \theta(g) = \theta(x_1^{r_1}\ldots x_n^{r_n}y_1^{s_1}\ldots y_m^{s_m})
= (x_1^{r_1}y_1^{s_1}\ldots y_m^{s_m}, \ldots, x_n^{r_n}y_1^{s_1}\ldots y_m^{s_m})$$
and so 
$$r_1 = \ldots = r_n = s_1 = \ldots = s_m = 0.$$
Thus $\theta$ is injective and so $G \in \la G_1,\ldots,G_n\ra$. Therefore we may assume that $G$ is defined by the presentation
$$\la x,y_1,\ldots,y_m \mid x^p = 1,\, y_j^{d_j} = 1, \, y_jy_{j'} = y_{j'}y_j,\, y_jxy_j\inv = x^{q_{j}}\, (j,j' \in [m])\ra.$$

In fact, we claim that we may assume that $d_j = d$ for every $j \in [m]$. Let $G'$ be defined by the presentation
\beq
\label{pac1}
\la x,y_1,\ldots,y_m \mid x^p = 1,\, y_j^{d} = 1, \, y_jy_{j'} = y_{j'}y_j,\, y_jxy_j\inv = x^{q_{j}}\, (j,j' \in [m])\ra.
\eeq
Note that $Q_{p,d_j} \subseteq Q_{p,d}$. Then $G$ is the image of $G'$ by the canonical homomorphism and so $G \in \la G'\ra$. 

Therefore we assume that $d_j = d$ for every $j \in [m]$, 
that is, $G$ is defined by the presentation (\ref{pac1}). Now let $r \in Q'_{p,d}$ and let $G_r$ be the group defined by the presentation obtained from (\ref{pac1}) by replacing the relation $y_1xy_1\inv = x^{q_1}$ by the relation $y_1xy_1\inv = x^{r}$. We show that
\beq
\label{pac2}
\mbox{there exists an injective homomorphism $\p:G \to G_r \times \Z/d\Z$.}
\eeq

If $q_1$ has order $d'$ in $\mathbb{F}_p^*$, then the subgroup generated by $r$ has a subgroup of order $d'$, which must necessarily contain $q_1$ because a cyclic group has a unique subgroup of a given order. Thus $q_1 \equiv r^k\,({\rm mod}\, p)$ for some $k \geq 0$. We define a homomorphism $\p:G \to G_r \times \Z/d\Z$ by $\p(x) = (x,0)$, $\p(y_1) = (y_1^k,1)$ and $\p(y_j) = (y_j,0)$ for $j = 2,\ldots,m$. We have 
$$\p(y_1xy_1\inv) = (y_1^kxy_1^{-k},0) = (x^{r^k},0) = (x^{q_1},0) = \p(x^{q_1})$$
and it follows easily that the homomorphism is well defined.

Let $g \in \ker(\p)$. Then $g = x^iy_1^{t_1}\ldots y_m^{t_m}$ for some $i \in \{0,\ldots,p-1\}$ and $t_1,\ldots,t_m \in \{ 0,\ldots,d-1\}$. Then
$$(1,0) = \p(g) = \p(x^iy_1^{t_1}\ldots y_m^{t_m}) = (x^iy_1^{kt_1}y_2^{t_2}\ldots y_m^{t_m},t_1),$$
hence $i = t_1 = \ldots = t_m = 0$ and $\p$ is injective. Therefore (\ref{pac2}) holds.

Now we assume that $G$ is defined by the presentation (\ref{pac1}) and 
we prove that $G \in \la G_{p,d}\ra$ by induction on $m$.
We prove first the case $m = 1$. Let $G$ be the group defined by the presentation
$$\la x,y \mid x^p = 1,\, y^{d} = 1, \, yxy\inv = x^q\ra.$$
This group is not necessarily isomorphic to $G_{p,d}$ because $q \in Q_{p,d}$ needs not be in $Q'_{p,d}$.
Let $r \in Q'_{p,d}$. By (\ref{pac2}), $G$ is isomorphic to a subgroup of $G_{p,d}^{(r)} \times \Z/d\Z = G_{p,d} \times \Z/d\Z$.
Since $\Z/d\Z \in \la G_{p,d}\ra$, it follows that $G \in \la G_{p,d}\ra$ and the case $m = 1$ is proved.

Assume now that $m > 1$ and the claim holds for $m-1$. In view of (\ref{pac2}), and exchanging $y_1$ and $y_2$ if needed, we may assume that $q_2 \in Q'_{p,d}$. Since $q_1 \in Q_{p,d}$ and cyclic groups have a unique subgroup of a given order, it follows that
there exists some $k \geq 0$ such that $q_1q_2^k \equiv 1\, ({\rm mod}\, p)$. Hence
$$y_2^ky_1xy_1\inv y_2^{-k} = y_2^kx^{q_1} y_2^{-k} = x^{q_1q_2^k} = x.$$
Now we replace $y_1$ by $y'_1 = y_2^ky_1$ in the generating set of $G$, adjusting the presentation accordingly. The unique relation that changes is $y_1xy_1\inv = x^{q_1}$, which is now $y'_1x(y'_1)\inv = x$. Thus $G$ admits the alternative presentation
$$\la x,y_1,y_2, \ldots,y_m \mid x^p = 1,\, y_j^{d} = 1, \, y_jy_{j'} = y_{j'}y_j,\, y_1xy_1\inv = x,\, y_rxy_r\inv = x^{q_{j}}\, (j,j',r \in [m],\, r > 1)\ra.$$
Denoting by $G'$ the group defined by the presentation
$$\la x,y_2,\ldots,y_m \mid x^p = 1,\, y_j^{d} = 1, \, y_jy_{j'} = y_{j'}y_j,\, y_jxy_j\inv = x^{q_{j}}\, (j,j' \in \{ 2,\ldots,m\})\ra,$$
it is now evident that $G \cong G' \times \Z/d\Z$. Since $G' \in \la G_{p,d}\ra$ by the induction hypothesis, we get $G\in \la G_{p,d}\ra$ and so $\ab(p) \ast \ab(d) = \la G_{p,d}\ra$.

Suppose now that $\ab(p) \ast \ab(d) = \la H\ra$ with $|H|$ minimum. The order of any element of any $G \in \ab(p) \ast \ab(d)$ must divide $|H|$: this follows from this property being preserved through (finitary) direct powers, subgroups and homomorphic images. Since $C_p, C_d \in \ab(p) \ast \ab(d)$, then $p$ and $d$ divide both $|H|$. Since $|G_{p,d}| = pd$, then $|H| = pd$ by minimality. 

Moreover, since $d < p$, it follows from the third Sylow theorem that $H$ has a normal subgroup of order $p$, generated, say, by $a$. Then $|H/\la a\ra| = d$, and, since $H \in \ab(p) \ast \ab(d)$, $H/\la a\ra$ is abelian.
Suppose that $H/\la a\ra$ is not cyclic. By the structure theorem of abelian groups, we have $H/\la a\ra \cong C_{d_1} \times \ldots \times C_{d_m}$ for some $d_1,\ldots,d_m > 1$ not all coprime. Hence $H/\la a\ra$ has exponent $d' < d$. For every $h \in H$, we have then $(h\la a\ra)^{d'} = \la a\ra$, yielding $h^{d'} \in \la a\ra$ and $h^{pd'} = 1$. Thus $H \models x^{pd'} = 1$ and so does $\ab(p) \ast \ab(d) = \la H\ra$. In particular, $G_{p,d} \models x^{pd'} = 1$. Using the presentation (\ref{presg}) of $G_{p,d}$, we obtain $b^{pd'} = 1$. Since $o(b) = d$ in $G_{p,d}$, then $d \mid pd'$ and therefore $d \mid d'$, contradicting $d' < d$. Thus $H/\la a\ra \cong C_d$.

Since $p$ and $d$ are coprime, it follows from the Schur-Zassenhaus theorem that $H \cong C_p \rtimes C_d$. Let $b \in H$ correspond to a generator of $C_d$. We have $a^p = b^d = 1$. Since $\la a\ra \lhd H$, we have $bab\inv = a^j$ for some $j \in [p-1]$. Applying Lemma \ref{diapr} to the case of $H$ (where $m = n = 1$), we get $j \in Q_{p,d}$, i.e.  
$j^d \equiv 1\,({\rm mod}\, p)$. It is now enough to prove that $o(j) = d$ in $\mathbb{F}_p^*$.

Let $d' = o(j)$ in $\mathbb{F}_p^*$. Since $j^d \equiv 1\,({\rm mod}\, p)$, we have $d' \mid d$. Suppose that $d' < d$. We show that  $H \models x^{d'}y^{d'} = y^{d'}x^{d'}$.

Indeed, $\la a\ra \lhd H$ implies that $h^{d'} \in \la a, b^{d'} \ra$ for every $h \in H$. Hence it suffices to show that $a$ and $b^{d'}$ commute in $H$. By (\ref{ndiapr1}), we have
$b^{d'}ab^{-d'} = a^{j^{d'}} = a$, hence $H \models x^{d'}y^{d'} = y^{d'}x^{d'}$ and so does $\ab(p) \ast \ab(d) = \la H\ra$. In particular, $G_{p,d} \models x^{d'}y^{d'} = y^{d'}x^{d'}$. Using the presentation (\ref{presg}) of $G_{p,d}$, we obtain $a^{d'}b^{d'} = b^{d'}a^{d'} = a^{d'q^{d'}}b^{d'}$. Hence $p \mid d'(q^{d'}-1)$. Since $d' < p$, it follows that $q^{d'} \equiv 1\,({\rm mod}\, p)$, contradicting $q \in Q'_{p,d}$. Therefore $o(j) = d$ in $\mathbb{F}_p^*$ as required.
\qed

Note that $|G_{p,d}| = pd$ constitutes a clear improvement with respect to $|F_2(p,d)| = p(d^2+1)d^2$.

\section{The pseudovariety $\U = \ds\bigvee_{p \in \mathbb{P}} \ab(p) \ast \ab(p-1)$}

This section is devoted to the pseudovariety $\U$, which is also the join of all the pseudovarieties $\ab(p) \ast \ab(d)$ studied in the preceding sections.

Let $\mathbf{E}$ denote the pseudovariety generated by all groups of prime order (equivalently, by all finite elementary abelian groups). 

Recall that {\bf Su} denotes the pseudovariety of finite supersolvable groups.

For every prime $p$, write $\U_p = \ab(p) \ast \ab(p-1)$. Hence $\U = \ds\bigvee_{p \in \mathbb{P}} \U_p$.

Given a pseudovariety $\V$ of finite groups, let
\bi
\item 
$\mathrm{Der}(\V)$ denote the class of all finite groups $G$ such that $[G,G] \in \V$;
\item 
$\mathrm{Syl}(\V)$ denote the class of all finite groups $G$ such that $H \in \V$ for every Sylow subgroup $H$ of $G$.
\ei

\bl
\label{dersyl}
Let $\V$ be a pseudovariety of finite groups. Then $\mathrm{Der}(\V)$ and $\mathrm{Syl}(\V)$ are pseudovarieties of finite groups.
\el

\proof
We must show that $\mathrm{Der}(\V)$ and $\mathrm{Syl}(\V)$ are closed under taking subgroups, homomorphic images and direct products.

Let $H \leq G \in \mathrm{Der}(\V)$. Then $[H,H] \leq [G,G] \in \V$ yields $H \in \mathrm{Der}(\V)$.

Assume now that $\p:G \to H$ is a surjective homomorphism and $G \in \mathrm{Der}(\V)$. Then $[G,G] \in \V$ yields 
$$[H,H] = [\p(G),\p(G)] = \p([G,G]) \in \V$$
and so $H \in \mathrm{Der}(\V)$.

Finally, assume that $G_1,G_2 \in \mathrm{Der}(\V)$. Then $[G_1,G_1], [G_2,G_2] \in \V$ yields 
$$[G_1 \times G_2,G_1 \times G_2] \leq [G_1,G_1] \times [G_2,G_2] \in \V$$
and consequently $G_1 \times G_2 \in \mathrm{Der}(\V)$.

Therefore $\mathrm{Der}(\V)$ is a pseudovariety of finite groups.

Now let $H \leq G \in \mathrm{Syl}(\V)$ and let $K$ be a Sylow subgroup of $H$. By the first Sylow theorem, $K$ is a subgroup of some Sylow subgroup $K'$ of $G$. Since $K' \in \V$, we get $K \in \V$, thus $H \in \mathrm{Syl}(\V)$.

Assume now that $\p:G \to H$ is a surjective homomorphism and $G \in \mathrm{Syl}(\V)$. Let $K$ be a Sylow subgroup of $H$. By \cite[Theorem 6.6]{Con2}, $K = \p(K')$ for some Sylow subgroup $K'$ of $G$. Now $K' \in \V$ yields $K \in \V$ and so $H \in \mathrm{Syl}(\V)$.

Finally, assume that $G_1,G_2 \in \mathrm{Syl}(\V)$ and let $K$ be a Sylow $p$-subgroup of $G_1 \times G_2$. Let $\pi_i:G_1 \times G_2 \to G_i$ be the canonical projection for $i = 1,2$. Then $\pi_1(K) \times \pi_2(K)$ is a $p$-subgroup of $G_1 \times G_2$ containing $K$, hence $K = \pi_1(K) \times \pi_2(K)$ by maximality of $K$. Again by the first Sylow theorem, $\pi_i(K)$ is a subgroup of some Sylow subgroup $K_i$ of $G_i$ for $i = 1,2$. It follows that $K_1, K_2 \in \V$, hence $\pi_1(K), \pi_2(K) \in \V$ and so $K = \pi_1(K) \times \pi_2(K) \in \V$. Thus $G_1 \times G_2 \in \mathrm{Syl}(\V)$.

Therefore $\mathrm{Syl}(\V)$ is also a pseudovariety of finite groups.
\qed

\bt
\label{charU}
The following conditions are equivalent for any finite group $G$:
\bi
\item[(i)] $G \in \U$;
\item[(ii)] $G \in {\bf Su} \cap \mathrm{Der}(\mathbf{E}) \cap \mathrm{Syl}(\ab)$;
\item[(iii)] $G$ embeds into $(C_{p_1} \rtimes C_{p_1-1}) \times \ldots \times (C_{p_s} \rtimes C_{p_s-1})$ for some $p_1,\ldots,p_s \in \mathbb{P}$, not necessarily distinct;
\item[(iv)] $G$ embeds into $(C_{p_1} \rtimes C_{p_1-1}) \times \ldots \times (C_{p_s} \rtimes C_{p_s-1}) \times C_{|G|}^m$ for some prime divisors $p_1,\ldots,p_s$ of $|G|$, not necessarily distinct, and some $m \in \N$.
\ei
\et

\proof
For every prime $p$, let $G_p$ denote the group $C_p \rtimes C_{p-1}$ presented by
$$\la x,y \mid x^p = 1,\, y^{p-1} = 1, \,  yxy\inv = x^{q_p}\ra,$$
for some fixed $q_p \in Q'_{p,p-1}$. 
By Proposition \ref{figen}, we have
$\U_p = \la G_{p}\ra$, hence
\beq
\label{charU1}
\U = \la G_{p} \mid p \in \mathbb{P} \ra.
\eeq
Write 
$$\U' = {\bf Su} \cap \mathrm{Der}(\mathbf{E}) \cap \mathrm{Syl}(\ab)$$
and let $\U''$ be the class of all finite groups satisfying condition (iii). 

(i) $\Rw$ (ii). 
We know that {\bf Su}, $\mathbf{E}$ and $\ab$ are pseudovarieties, and so are $\mathrm{Der}(\mathbf{E})$ and $\mathrm{Syl}(\ab)$ in view of Lemma \ref{dersyl}. Hence $\U'$ is a pseudovariety. In view of (\ref{charU1}), it suffices to show that $G_p \in \U'$ for every $p \in \mathbb{P}$.

Let $p \in \mathbb{P}$. We have $[G_p,G_p] \leq \la x \ra \cong C_p$, hence $G_p \in \mathrm{Der}(\mathbf{E})$. On the other hand, $G_p \in {\bf G}_p \ast \ab(p-1)$ yields $G_p \in {\bf Su}$ by \cite[Corollary 2.7]{AS}. 

Finally, let $H \leq G_p$ be a $q$-Sylow subgroup for some $q \in \mathbb{P}$. If $q = p$, then $H = \la x \ra \in \ab$, hence we may assume that $q \neq p = o(x)$. Then
$$[H,H] \leq [G_p,G_p] \cap H \leq \la x\ra \cap H = \{1\},$$
hence $H \in \ab$ and so $G_p \in \mathrm{Syl}(\ab)$. Therefore $G_p \in \U'$ as required.

(ii) $\Rw$ (iii). By Dirichlet's theorem, every finite cyclic group embeds into $C_{p-1}$ for some $p \in \mathbb{P}$. Now it follows from the fundamental theorem of finite abelian groups that $\ab \subseteq \U''$.

Suppose that $\U' \not\subseteq \U''$ and choose $G \in \U'\setminus \U''$ with $|G|$ minimal. We show that
\beq
\label{charU2}
\mbox{if $1 \neq N,N' \unlhd G$, then $N \cap N' \neq 1$.} 
\eeq
Indeed, suppose that $N \cap N' = 1$. Then the mapping
$$\begin{array}{rcl}
\p:G&\to&G/N \times G/N'\\
g&\to&(Ng,N'g)
\end{array}$$
is injective. Since $G \in \U'$, we have $G/N, G/N' \in \U'$. Since $|G/N|, |G/N'| < |G|$, it follows from the minimality of $G$ that $G/N,G/N' \in \U''$. But then $G$ would be in $\U''$, since $\p$ is an embedding into a direct product of groups in $\U''$. We have reached a contradiction, therefore (\ref{charU2}) holds.

For each $p \in \mathbb{P}$, let $O_p(G)$ denote the intersection of all the Sylow $p$-subgroups of $G$. This is a normal subgroup of $G$, which is in fact nilpotent (being a $p$-group). The {\em Fitting subgroup} $F(G)$ is the (unique) largest normal nilpotent subgroup of $G$. Hence $O_p(G) \leq F(G)$ for every $p \in P$.

By the first Sylow theorem, $O_p(G)$ is a subgroup of some Sylow subgroup $S_p$ of $F(G)$, which is itself a subgroup of some Sylow subgroup $S'_p$ of $G$. In view of the second Sylow theorem, we get
 $$O_p(G) = \bigcap_{g\in G} gO_p(G)g\inv  \leq \bigcap_{g\in G} gS_pg\inv \leq \bigcap_{g\in G} gS'_pg\inv = O_p(G),$$
 hence $\bigcap_{g\in G} gS_pg\inv = O_p(G)$. On the other hand, since $F(G) \unlhd G$, we get $gS_pg\inv \leq F(G)$ for every $g \in G$, and since a nilpotent group has a unique Sylow $p$-subgroup, $gS_pg\inv = S_p$ for every $g \in G$. Thus $O_p(G)$ is the (unique) Sylow $p$-subgroup of $F(G)$. Since a nilpotent group is always the direct product of its Sylow subgroups, we get 
$$F(G) = \prod_{p\in\mathbb{P}} O_p(G).$$
Since $[G,G]$ is abelian, and therefore nilpotent, we have $[G,G] \leq F(G)$. Since $\ab \subseteq \U''$, then $G$ is not abelian and so $F(G) \neq 1$. By (\ref{charU2}), there is at most one prime $p$ such that $O_p(G)\neq 1$, thus $F(G) = O_p(G) = S_p$ for this only prime $p$, which we now fix.

Since $G\in {\rm Syl}(\mathbf{Ab})$ then $F(G)=O_p(G)=S_p$ is abelian. Since $G$ is a finite solvable group,  $C_G(F(G))\leq F(G)$ (see \cite[Theorem 6.1.3]{Gor}). We therefore obtain $C_G(O_p(G))=O_p(G)$.

Let $Y = \Omega_1(O_p(G))$ (i.e. the group generated by all elements of $O_p(G)$ of order $p$). Since $O_p(G)$ is abelian, then $Y$ is elementary abelian. Moreover, $G$ acts on $Y$ by conjugation.
It follows that $Y$ is a $G$-module over the finite field $\mathbb{F}_p$, and the $G$-submodules of $Y$ are exactly the $G$-normal subgroups of $Y$. Since 
$C_G(Y) \geq O_p(G)$, then 
$G/C_G(Y)$ is a quotient of $G/O_p(G)$, hence $G/C_G(Y)$ has order coprime to $p$ (since $O_p(G)$ is a Sylow $p$-subgroup of $G$). Now it follows from Maschke's theorem that $Y$ can be written in the form $Y=Y_1 \times \ldots \times Y_e$, where each $Y_i$ is an irreducible $G$-submodule of $Y$ (that is, a minimal $G$-normal subgroup of $Y$).

Thus each $Y_i$ is a chief factor of $G$. Since $G$ is supersolvable, each $Y_i$ has prime order, hence $Y_i$ is cyclic of order $p$. If $e \geq 2$, then it would follow that $Y_1$ and $Y_2$ are two normal subgroups of $G$ intersecting trivially, contradicting (\ref{charU2}). Thus $Y$ is cyclic of order $p$. Since $O_p(G)$ is abelian, it follows that $O_p(G)$ is itself cyclic.

Since $[G,G] \leq F(G) = O_p(G)$, the quotient group $A = G/O_p(G)$ is abelian. 
Since $|A|$ is coprime to $p$, it follows from Schur-Zassenhaus theorem that we can write $G = O_p(G) \rtimes A$. Thus there exists some $A' \leq G$, isomorphic to $A$, such that $G = O_p(G)A'$, and the action of $A'$ on $O_p(G)$ is by conjugation.

Denoting by $C_{A'}(O_p(G)) = \{ a \in A'\mid g^a = g \, \forall g \in O_p(G)\}$ the centralizer of this action, and since $C_G(O_p(G))=O_p(G)$, we have
$$C_{A'}(O_p(G)) = A' \cap C_G(O_p(G)) = A' \cap O_p(G) = 1,$$  
hence the action of $A$ on $O_p(G)$ is faithful and so
$A$ embeds in the automorphism group of $O_p(G)$.

Since $A'$ acts by automorphisms (in fact, by conjugation) on the abelian group $O_p(G)$ and $|A'|$ is coprime to $|O_p(G)|$, it follows from a theorem of Fitting \cite[Theorem 4.34]{Isa} that we also have $O_p(G) = [O_p(G),A'] \times C_{A'}(O_p(G))$, hence 
$$O_p(G) = [O_p(G),A'] \times C_{A'}(O_p(G)) =  [O_p(G),A'] \leq [G,G] \leq F(G) = O_p(G),$$
yielding $O_p(G) = [G,G]$. Since $O_p(G)$ is a cyclic $p$-group and $[G,G] \in {\bf E}$, it follows that $O_p(G) \cong C_{p}$. Since $A$ embeds in the automorphism group of $O_p(G)$, then $A$ embeds into $C_{p-1}$.

But then $G = O_p(G) \rtimes A \in \U''$, a contradiction. Therefore $\U' \subseteq \U''$ and (iii) holds.

(iii) $\Rw$ (iv). Assume there exists an embedding $\p:G \to G_{p_1} \times \ldots \times G_{p_s}$ for some $p_1,\ldots,p_s \in \mathbb{P}$. For $i \in [s]$, let $\pi_i: G_{p_1} \times \ldots \times G_{p_s} \to G_{p_i}$ be the canonical projection. Then
\beq
\label{charU3}
\p(G) \leq \pi_1(\p(G)) \times \ldots \times \pi_s(\p(G)),
\eeq
and we may assume that $\pi_i(\p(G)) \neq 1$ for every $i \in [s]$.

Suppose that $p_i$ does not divide $|G|$ for some $i$. We claim that $\pi_i(\p(G)) \in \ab$. Indeed, if $\pi_i(\p(G)) \not\in \ab$, then $\pi_i(\p(G))$ contains a nontrivial commutator $c$ of $G_{p_i}$, which has necessarily order $p_i$. Let $g \in G$ be such that $\pi_i(\p(g)) = c$. Then
$$c^{|G|} = (\pi_i(\p(g)))^{|G|} = \pi_i(\p(g^{|G|})) = 1$$
and so $p_i$ divides $|G|$, a contradiction. Thus $\pi_i(\p(G)) \in \ab$. Since $|\pi_i(\p(G))|$ divides $|G|$, it follows from the fundamental theorem of finite abelian groups that $\pi_i(\p(G))$ embeds into $C_{|G|}^{m_i}$ for some $m_i \in \N$.
In view of (\ref{charU3}), we get the desired embedding.

(iv) $\Rw$ (i). By Dirichlet's theorem, $C_{|G|}$ embeds into $C_{p-1}$ for some $p \in \mathbb{P}$ and we are done.
\qed

We immediately get:

\bc
\label{mud}
The pseudovariety $\U$ has decidable membership problem.
\ec

We recall that $\M$ denotes the pseudovariety of all finite metabelian groups.

\bp
\label{UinM}
$\ab \subset \U \subset {\bf Su} \cap \M$.
\ep

\proof
Since $\M = \mathrm{Der}(\ab) \supseteq \mathrm{Der}({\bf E})$, Theorem \ref{charU} yields $\U \subseteq {\bf Su} \cap \M$.

Let $Q_8$ denote the quaternion group, which is well known to be both supersolvable and metabelian. Since $Q_8$ is not abelian and is its unique Sylow subgroup, it follows from Theorem \ref{charU} that $Q_8 \notin \U$. Therefore 
$\U \subset {\bf Su} \cap \M$.

The inclusion $\ab \subset \U$ follows from Theorem \ref{charU} and is obviously strict.
\qed

\section{The variety generated by $\U$}

The main result of this section proves that $\U$ generates the variety $\cal{M}$ of all metabelian groups. It is known that
$\cal{M}$ is generated by $F_2({\cal{M}})$, the free metabelian group of rank 2 (see the comments before 16.41 in \cite{Neu}). 

On the other hand, it follows from (\ref{charU1}) that 
\beq
\label{vugp}
{\cal{V}}(\U) = {\cal{V}}(\{ G_p \mid p \in \mathbb{P} \}).
\eeq

We present now a concrete description of $F_2({\cal{M}})$. We fix an alphabet $A = \{ a,b\}$. The Cayley graph Cay$_A(\mathbb{Z} \times \mathbb{Z})$ is an infinite grid with vertex set $\mathbb{Z} \times \mathbb{Z}$ and edges of the form
$$(m,n) \xr{a} (m+1,n),\quad (m,n) \xr{b} (m,n+1)\quad (m,n \in \mathbb{Z}).$$
We denote by $E$ the set of all such edges. Each edge can be read in both directions: in the direct sense, with label $a$ or $b$; in the opposite sense, with label $a\inv$ or $b\inv$ (so we may consider edges labelled by $a\inv$ or $b\inv$ whenever convenient). We consider vertex $(0,0)$ as the basepoint of Cay$_A(\mathbb{Z} \times \mathbb{Z})$. It is immediate that each word $u \in (A \cup A\inv)^*$ labels a unique path $\pi_u$ off the basepoint in Cay$_A(\mathbb{Z} \times \mathbb{Z})$. For every $e \in E$, we denote by $|u|_e$ the number of times $\pi_u$ traverses $e$ in the positive direction minus the number of times $\pi_u$ traverses $e$ in the opposite direction. Note that if $\oo{u}$ denotes the reduced form of $u$, then $|u|_e = |\oo{u}|_e$ and so $|g|_e$ is well defined for every $g \in F_2$. It follows from a general result of Almeida \cite[Theorem 2.2]{Alm} (see also \cite{LS}) that $F_2({\cal{M}}) = F_2/N$ for
$$N = \{ g \in F_2 \;\big{|}\; |g|_e = 0 \mbox{ for every }e \in E\}.$$

Given a prime $p \in \mathbb{P}$, we say that $q \in \mathbb{Z}$ is a {\em primitive root modulo} $p$ if $q$ corresponds to a generator of $\mathbb{F}_p^*$. This is equivalent to saying that the remainder of the division of $q$ by $p$ is in $Q'_{p,p-1}$.
Given $q \in \mathbb{Z}$, we write
$${\rm PR}(q) = \{ p \in \mathbb{P} \mid q \mbox{ is a primitive root modulo }p\}.$$

\bt
\label{vum}
We have ${\cal{V}}(\U) = \cal{M}$.
\et

\proof
It follows from Proposition \ref{UinM} that $\U \subset \M \subset \cal{M}$, hence ${\cal{V}}(\U) \subseteq \cal{M}$.

To prove the converse inclusion, since ${\cal{M}} = {\cal{V}}(F_2({\cal{M}}))$, it suffices to show that $F_2({\cal{M}}) = F_2/N  \in {\cal{V}}(\U)$, and this will follow from:
\beq
\label{vum1}
\begin{array}{r}
\mbox{for every $u \in F_2 \setminus N$, there exists some $p_u \in \mathbb{P}$ and some homomorphism}\\
 \p_u:F_2/N \to G_{p_u}\mbox{ such that $\p_u(uN) \neq 1$.}
 \end{array}
\eeq
Indeed, if (\ref{vum1}) holds, then 
$$\p:F_2/N \to \prod_{u \in F_2 \setminus N} G_{p_u}$$
defined by $\p(wN) = (\p_u(wN))$ is an injective homomorphism. Since this direct product is in ${\cal{V}}(\U)$ by (\ref{vugp}), the theorem follows.

To prove (\ref{vum1}), we start by simplifying the problem. We say that $u$ is {\em first-quadrant} if all the vertices occurring in the path $\pi_{\oo{u}}$ have both coordinates non-negative. We claim that it is enough to show that (\ref{vum1}) holds for first-quadrant words. Indeed, assume that this is true, and let $u \in F_2 \setminus N$ be arbitrary. Then there exists some $m \geq 0$ such that $v = a^mb^mub^{-m}a^{-m}$ is first-quadrant (and $v \in F_2 \setminus N$ since $N \lhd F_2$). If $\p_v:F_2/N \to G_{p_v}$ satisfies $\p_v(vN) \neq 1$, then also $\p_v(uN) \neq 1$, therefore it suffices to prove (\ref{vum1}) for first-quadrant words. 

For every $n \in \mathbb{Z}$, let
$h_n(u) = \ds\sum_{m \in \mathbb{Z}} |u|_{e_m}$, where $e_m$ is the edge $(m,n) \xr{a} (m+1,n)$. Let also $v_n(u) = \ds\sum_{m \in \mathbb{Z}} |u|_{e'_m}$, where $e'_m$ is the edge $(n,m) \xr{b} (n,m+1)$. 
Next, we prove (\ref{vum1}) for the particular case where $h_n(u) \neq 0$ for some $n$. We shall assume that $n$ is maximum. Note that $n \geq 0$ since $u$ is first-quadrant.

Let $k = |u|$. By \cite[Corollary 2]{HB}, there are at most two positive primes $p$ such that PR$(p)$ is finite. Thus there exists some prime $q > k$ such that PR$(q)$ is infinite. In particular, we can fix some $p \in {\rm PR}(q)$ such that $p > kq^k$. 

Since $q \in Q'_{p,p-1}$, then $G_p$ can be presented by
$$\la a,b \mid a^p = 1,\, b^{p-1} = 1, \,  bab\inv = a^{q}\ra.$$
Since $G_p \in \cal{M}$ and $F_2/N$ is the free metabelian group on $A$, we have a homomorphism $\p:F_2/N \to G_p$ extending the identity map on $A$. We claim that $\p(uN) \neq 1$.

Indeed, let $u_i$ denote the prefix of length $i$ of $u$ for $i = 0,\ldots,k$. Suppose that $(m_0,n_0)$ is the end vertex of $\pi_{u_i}$, and let
\beq
\label{vum2}
(m_1,n_1) \xr{a^{\varepsilon_1}} (m_1+\varepsilon_1,n_1), \ldots, (m_s,n_s) \xr{a^{\varepsilon_s}} (m_s+\varepsilon_s,n_s)
\eeq
be all the edges of $\pi_{u_i}$ labelled by either $a$ or $a\inv$. We prove that 
\beq
\label{newvum1}
u_i = a^{\sum_{j=1}^s \varepsilon_jq^{n_j}}b^{n_0}
\eeq
holds in $G_p$ for $i= 0,\ldots,k$ by induction.

The claim holds trivially for $i= 0$. Let $i \in \{ 0,\ldots,k-1\}$ and assume the claim holds for $i$. If $u_{i+1} = u_ib^{\varepsilon}$ for $\varepsilon = \pm 1$, then we have no more edges to add to (\ref{vum2}) going from $\pi_{u_i}$ to 
$\pi_{u_{i+1}}$, but the end vertex changes to $(m_0,n_0+\varepsilon)$. Using the induction hypothesis, the equalities
$$u_{i+1} = u_ib^{\varepsilon} = a^{\sum_{j=1}^s \varepsilon_jq^{n_j}}b^{n_0}b^{\varepsilon} = a^{\sum_{j=1}^s \varepsilon_jq^{n_j}}b^{n_0+\varepsilon}$$
hold in $G_p$ as required.

Thus we may assume that $u_{i+1} = u_ia^{\varepsilon}$ for $\varepsilon = \pm 1$. Comparing $\pi_{u_{i+1}}$ with $\pi_{u_{i}}$, we must add the edge $(m_0,n_0) \xr{a^{\varepsilon}} (m_0+\varepsilon,n_0)$ to (\ref{vum2}), and $(m_0+\varepsilon,n_0)$ is now the new end vertex. Since $u$ is first-quadrant, we have $n_0 \geq 0$ and it follows from previous computations that $b^{n_0}ab^{-n_0} = a^{q^{n_0}}$ holds in $G_p$.
Using the induction hypothesis, we now get the equalities
$$u_{i+1} = u_ia^{\varepsilon} = a^{\sum_{j=1}^s \varepsilon_jq^{n_j}}b^{n_0}a^{\varepsilon} = a^{\sum_{j=1}^s \varepsilon_jq^{n_j}+\varepsilon q^{n_0}}b^{n_0}$$
to hold in $G_p$ and (\ref{newvum1}) is verified.

In particular, (\ref{newvum1}) holds for $u_k = u$. Recall that $h_n(u) \neq 0$ for maximum $n$. 
Assume that
$$(m_1,n) \xr{a^{\varepsilon_1}} (m_1+\varepsilon_1,n), \ldots, (m_s,n) \xr{a^{\varepsilon_s}} (m_s+\varepsilon_s,n),$$

$$(m'_1,n_1) \xr{a^{\varepsilon'_1}} (m'_1+\varepsilon'_1,n_1), \ldots, (m'_t,n_t) \xr{a^{\varepsilon'_t}} (m'_t+\varepsilon'_t,n_t)$$
and
$$(m''_1,n'_1) \xr{a^{\varepsilon''_1}} (m''_1+\varepsilon''_1,n'_1), \ldots, (m''_r,n'_r) \xr{a^{\varepsilon''_r}} (m''_r +\varepsilon''_r,n'_r)$$
are now all the edges of $\pi_{u}$ labelled by either $a$ or $a\inv$, with $n_1,\ldots,n_t < n$ and $n'_1,\ldots,n'_r > n$. By our preceding claim, the equality 
$$u =  a^{\sum_{j=1}^s \varepsilon_jq^{n}+\sum_{\ell =1}^t \varepsilon'_{\ell}q^{n_{\ell}}+\sum_{h =1}^r \varepsilon''_{h}q^{n'_{h}}}b^{n_0}$$
holds in $G_p$ for some $n_0$. Let $S = \sum_{j=1}^s \varepsilon_jq^{n}$, $T = \sum_{\ell =1}^t \varepsilon'_{\ell}q^{n_{\ell}}$ and $R = \sum_{h =1}^r \varepsilon''_hq^{n'_{h}}$. By maximality of $n$, we have $R = 0$.
Now on the one hand, we have
$$|T| \leq \sum_{\ell =1}^t q^{n_{\ell}} \leq tq^{n-1} \leq kq^{n-1}.$$
On the other hand, we have $S = h_n(u)q^n$ and so $h_n(u) \neq 0$ yields $q^n \leq |S| \leq kq^n$. Hence
$$0 < q^n-kq^{n-1} \leq |S|-|T| \leq |S+T| \leq kq^n \leq kq^k < p.$$
Since $a$ has order $p$ in $G_p$, it follows from (\ref{newvum1})
that $a^{R+S+T}b^{n_0} = a^{S+T}b^{n_0} \neq 1$ in $G_p$, that is, $u \neq 1$ in $G_p$. Therefore $\p(uN) \neq 1$.

Out of symmetry, (\ref{vum1}) also holds when $v_n(u) \neq 0$ for some $n$.

Finally, we prove (\ref{vum1}) in full generality. Let $u \in F_2 \setminus N$. Then $|u|_e \neq 0$ for some $e \in E$. Out of symmetry, we may assume that $e$ is an edge of the form $(m,n) \xr{a} (m+1,n)$, with $n$ maximum and $m$ maximum for such $n$. Recall that $u$ is first-quadrant, so $m,n \geq 0$. Our goal is to find some endomorphism $\theta$ of $F_2/N$ such that $\theta(uN) = vN$ for some first-quadrant $v \in F_2$ satisfying $h_c(v) \neq 0$ for some $c \geq 0$. By the preceding case, we know that there exists some $p \in \mathbb{P}$ and some homomorphism $\p:F_2/N \to G_{p}$ such that $\p(vN) \neq 1$. Then the homomorphism $\p\circ \theta:F_2/N \to G_{p}$ satisfies $(\p\circ\theta)(uN) = \p(vN) \neq 1$ and we are done. 

So let $k = |u|$. Let $\theta_0$ be the endomorphism of $F_2$ defined by $\theta_0(a) = ab$ and $\theta_0(b) = b^k$. Since $u$ is first-quadrant, so is $v = \theta_0(u)$. Since $F_2/N$ is the free metabelian group over $A = \{ a,b\}$, then $\theta_0$ induces a well-defined endomorphism $\theta$ of $F_2/N$. How does $\theta_0$ act on $E$? It is easy to check that each occurrence of an edge $(i,j) \xr{a} (i+1,j)$ in $\pi_u$ induces the occurrence of a subpath 
$$(i,i+jk) \xr{a} (i+1,i+jk) \xr{b} (i+1,i+1+jk)$$
in $\pi_v$. On the other hand, each occurrence of an edge $(i,j) \xr{b} (i,j+1)$ in $\pi_u$ induces the occurrence of a subpath $(i,i+jk) \xr{b^k} (i,i+(j+1)k)$ in $\pi_v$. 

We claim that $h_{m+nk}(v) = |u|_e \neq 0$. Indeed, each occurrence of $e$ in $\pi_u$ contributes an occurrence of $(m,m+nk) \xr{a} (m+1,m+nk)$, in the same direction. Thus it suffices to show that no other edge occurring in $\pi_u$ contributes to the value of $h_{m+nk}(v)$. We only have to consider edges labelled by $a$. Let $f$ be the edge $(i,j) \xr{a} (i+1,j)$. If $j > n$, it follows from the maximality of $n$ that $|u|_{f} = 0$, so the occurrences of $(i,i+jk) \xr{a} (i+1,i+jk)$ in $\pi_v$ would match each other. The same would happen if $f$ is the edge $(i,n) \xr{a} (i+1,n)$ with $i > m$, by maximality of $m$.

Suppose now that $f$ is an edge $(i,n) \xr{a} (i+1,n)$ with $i < m$. Then the induced occurrences in $\pi_v$ feature the edge $(i,i+nk) \xr{a} (i+1,i+nk)$. Since $i+nk < m+nk$, we have no contributions for the value of $h_{m+nk}(v)$. Finally, if 
$f$ is an edge $(i,j) \xr{a} (i+1,j)$ with $j < n$, then the induced occurrences in $\pi_v$ feature the edge $(i,i+jk) \xr{a} (i+1,i+jk)$. Since $i+jk < k + (n-1)k \leq m+nk$, we have no contributions for the value of $h_{m+nk}(v)$ either. Thus  
$h_{m+nk}(v) = |u|_e \neq 0$ as required.
\qed

Recall that, if $\cal{V}$ is a variety of groups, we denote by ${\cal{V}}^f$ the pseudovariety of all finite groups in $\cal{V}$. Such pseudovarieties have some positive features (see \cite{MSTm}), but this is not the case of $\U$:

\bc
\label{novfin}
There exists no variety $\cal{V}$ of groups satisfying $\U = {\cal{V}}^f$. 
\ec

\proof
Indeed, suppose that $\U = {\cal{W}}^f$ for some variety $\cal{W}$ of groups. Then $\U \subset \cal{W}$ and Theorem \ref{vum} yields ${\cal{M}} = {\cal{V}}(\U) \subseteq \cal{W}$. Thus $\M = {\cal{M}}^f \subseteq {\cal{W}}^f = \U$, contradicting Proposition \ref{UinM}.
\qed

Given $q \in \mathbb{P}$, the {\em Baumslag-Solitar group} $BS(1,q)$ is defined by the presentation
$$\langle a,b \mid bab\inv = a^q\rangle.$$
Baumslag-Solitar groups were introduced in \cite{BSo} and have deserved full attention in geometric group theory ever since, constituting often the source of interesting counter-examples.

We can prove the following:

\bt
\label{bsq}
The following assertions hold:
\bi
\item[(i)]
The equality 
\beq
\label{bsq1}
{\cal{V}}(BS(1,q)) = {\cal{V}}(\{ G_p \mid p \in {\rm PR}(q) \}).
\eeq
fails for at most two primes $q \in \mathbb{P}$. 
\item[(ii)]
The equality (\ref{bsq1}) holds for $q \in \mathbb{P}$ if the Riemann hypothesis holds for the Dedekind zeta function of each field $\mathbb{Q}(\sqrt[p]{q})$ for every $p \in \mathbb{P}$.
\ei
\et

\proof
We show that (\ref{bsq1}) holds whenever ${\rm PR}(q)$ is infinite. Then part (i) follows from \cite[Corollary 2]{HB}. On the other hand, part (ii) follows then from a theorem of Hooley \cite{Hoo}, who showed that ${\rm PR}(q)$ is infinite if the Riemann hypothesis holds for the Dedekind zeta function of each field $\mathbb{Q}(\sqrt[p]{q})$ for every $p \in \mathbb{P}$.

Assume then that ${\rm PR}(q)$ is infinite. Let $p \in {\rm PR}(q)$. Then
$$\langle a,b \mid a^p = b^{p-1} = 1, \, bab\inv = a^q\rangle$$
is a presentation for $G_p$ and so we have a canonical homomorphism from $BS(1,q)$ onto $G_p$. Hence $G_p \in {\cal{V}}(BS(1,q))$ for every $p \in {\rm PR}(q)$ and so ${\cal{V}}(\{ G_p \mid p \in {\rm PR}(q) \}) \subseteq {\cal{V}}(BS(1,q))$.

For the converse, let $(p_1,p_2,\ldots)$ be an enumeration of the elements of ${\rm PR}(q)$. For each $n \geq 1$, let $\p_n:BS(1,q) \to G_{p_n}$ be the canonical homomorphism. Then
$$\begin{array}{rcl}
\p:BS(1,q)&\to&\ds\prod_{n\geq 1} G_{p_n}\\
g&\mapsto&(\p_n(g))
\end{array}$$
is a homomorphism.
We claim that $\p$ is injective. Let $g \in \ker(\p)$. 
Since $BS(1,q)$ is an HNN extension, it follows from Britton's normal form (see \cite{LSc}) that $g = a^ib^ja^{k}$ for some $i,j,k \in \mathbb{Z}$. Since $g \in \ker(\p)$ if and only if $g\inv \in \ker(\p)$, we may assume that $j \geq 0$. Since ${\rm PR}(q)$ is infinite, there exists some $n \geq 1$ such that $p_n > \max(|i+kq^j|,j+1)$. Then
$$1 = \p_n(g) = \p_n(a^ib^ja^{k}) = a^ib^ja^{k} = a^{i+kq^j}b^j.$$
Since $o(a) = p_n$ and $o(b) = p_n-1$ in $G_{p_n}$, it follows from our choice of $p_n$ that $i+kq^j = j = 0$. Hence $k = -i$ and so $g =  a^ib^ja^{k} = a^{i}a^{-i} = 1$. Thus $\p$ is injective and so ${\cal{V}}(BS(1,q)) \subseteq {\cal{V}}(\{ G_p \mid p \in {\rm PR}(q) \})$. Therefore ${\cal{V}}(BS(1,q)) = {\cal{V}}(\{ G_p \mid p \in {\rm PR}(q) \})$ as required.
\qed

Given $P \subseteq \mathbb{P}$, let
$${\cal{M}}_P = {\cal{V}}(\{ G_p \mid p \in {\rm PR}(q) \}).$$

It would be tempting to think that ${\cal{M}}_P \neq {\cal{M}}_Q$ for different subsets $P,Q \subseteq \mathbb{P}$, but that is far from true. Indeed, it follows from the theorems of Cohen \cite{Coh} that there exist only countably many varieties of metabelian groups. Hence there exists an uncountable subset $\cal{P}$ of $2^{\mathbb{P}}$ such that ${\cal{M}}_P = {\cal{M}}_Q$ for all $P,Q \in \cal{P}$. Since $\cal{P}$ is uncountable, there must exist some $P,Q \in {\cal{P}}$ with $P\setminus Q$ infinite. Let $R = \mathbb{P} \setminus (P\cup Q)$. Then
$${\cal{M}} = {\cal{M}}_{\mathbb{P}} = {\cal{M}}_{P \cup Q \cup R} = {\cal{M}}_{Q \cup R}.$$
Since $\mathbb{P} \setminus (Q\cup R) = P\setminus Q$, it follows that there exists some $S \subset \mathbb{P}$ such that ${\cal{M}} = {\cal{M}}_{S}$ and $\mathbb{P} \setminus S$ is infinite.

This suggests then the following conjecture:

\bcon
\label{conjecture}
If $P \subseteq \mathbb{P}$ is infinite, then ${\cal{M}}_P = {\cal{M}}$.
\econ

If the conjecture holds, it follows from Theorem \ref{bsq}(i) that the equality ${\cal{V}}(BS(1,q)) = {\cal{M}}$ fails for at most two primes $q$.

\section{The pro-$\U$ topology}

We discuss in this section the pro-$\U$ topology, starting with
a general result concerning the $\U$-closure of a finitely generated subgroup. However, it does not provide decidability for being finitely generated, neither computability in that case.

\bp
\label{ucl}
Let $H \leq_{f.g.} F_n$. Then
$$\mathrm{Cl}_{\U}(H) = \bigcap_{p\in \mathbb{P}} \mathrm{Cl}_{\U_p}(H) = \bigcap_{p\in \mathbb{P}} HL_{n,p,p-1} \geq H\mathrm{Cl}_{\U}(1) = H\bigcap_{p\in \mathbb{P}} L_{n,p,p-1} \geq HF_n^{(2)} = \mathrm{Cl}_{\M}(H).$$
\ep

\proof
Since $\U_p \subseteq \U$, it follows from (\ref{icl}) and the fact that the closure of a subgroup is always a subgroup that $\mathrm{Cl}_{\U}(H) \leq \mathrm{Cl}_{\U_p}(H)$ for every $p \in \mathbb{P}$, therefore
$$\mathrm{Cl}_{\U}(H) \leq \bigcap_{p\in \mathbb{P}} \mathrm{Cl}_{\U_p}(H).$$
We prove next the opposite inclusion. Write $F = F_n$ and $C = \bigcap_{p\in \mathbb{P}} \mathrm{Cl}_{\U_p}(H)$. It is well known (see e.g. \cite[Proposition 1.3]{MSW}) that we have
$$\mathrm{Cl}_{\U}(H) = \bigcap\{ N \unlhd F \mid F/N \in \U \mbox{ and }H \leq N\}.$$
Hence it suffices to show that $C \leq N$ for every normal  subgroup $N$ of $F$ containing $H$ with $F/N\in \mathbf{U}$.

Now, it follows from Theorem \ref{charU} that $F/N$ embeds into $G_{p_1} \times \ldots \times G_{p_s}$ for some $p_1,\ldots,p_s \in \mathbb{P}$. For $i \in [s]$, let $\pi_i:G_{p_1} \times \ldots \times G_{p_s} \to G_{p_i}$ be the canonical projection, and let $\p:F \to F/N$ be the canonical homomorphism. Let $N_i = \ker(\pi_i\circ \p)$. We claim that
\beq
\label{ucl1}
\ds\bigcap_{i=1}^s N_i \leq N.
\eeq 
Let $u \in \ds\bigcap_{i=1}^s N_i$. Then $\pi_i(\p(u)) = 1$ for every $i \in [s]$, hence $\p(u) = 1$ and so $u \in N$. Thus (\ref{ucl1}) holds.

For every $i \in [s]$, we have
$$F/N_i = F/\ker(\pi_i\circ \p) \cong \mathrm{Im}(\pi_i\circ \p) \leq G_{p_i} \in \U_{p_i},$$
hence $F/N_i \in \U_{p_i}$. On the other hand, given $i \in [s]$, $H \leq N$ yields
$\pi_i(\p(H)) = \pi_i(1) = 1$, thus $H \leq N_i$. It follows that $\mathrm{Cl}_{\U_{p_i}}(H) \leq N_i$ for every $i \in [s]$, hence (\ref{ucl1}) yields
$$C = \bigcap_{p\in \mathbb{P}} \mathrm{Cl}_{\U_p}(H) \leq \bigcap_{i=1}^s \mathrm{Cl}_{\U_{p_i}}(H)  \leq \ds\bigcap_{i=1}^s N_i \leq N$$
and so the equality $\mathrm{Cl}_{\U}(H) = \bigcap_{p\in \mathbb{P}} \mathrm{Cl}_{\U_p}(H)$ holds.

The second and third equalities follow from Proposition \ref{clop}. The fourth equality follows from \cite[Theorem 4.7(i)]{MSTm}.

Now $H\mathrm{Cl}_{\U}(1) \leq \mathrm{Cl}_{\U}(H)$ follows from $H \leq \mathrm{Cl}_{\U}(H)$ and $\mathrm{Cl}_{\U}(1) \leq \mathrm{Cl}_{\U}(H)$. 

Finally, $\U \subseteq \M$ yields $\mathrm{Cl}_{\M}(1) \leq \mathrm{Cl}_{\U}(1)$ and consequently $F^{(2)} \leq \mathrm{Cl}_{\U}(1)$ by \cite[Theorem 4.7(i)]{MSTm}. Thus $HF^{(2)} \leq H\mathrm{Cl}_{\U}(1)$ and the proof is complete.
\qed

We can prove however, for an arbitrary pseudovariety and any finite index subgroup of an arbitrary group, that the closure is computable. But we need a lemma first.

\bl
\label{fi}
{\rm \cite[Proposition 2.6]{MSTm}} 
Let $\V$ be a pseudovariety of finite groups. Let $G$ be a group and let $H \leq G$ have finite index.
Then the following conditions are equivalent:
\bi
\item[(i)]
$H$ is $\V$-closed;
\item[(ii)]
$H$ is $\V$-clopen;
\item[(iii)]
$G/{\rm Core}_{G}(H) \in \V$.
\ei
\el

\bp
\label{fid}
Let $\V$ be a pseudovariety of finite groups with decidable membership problem. Let $G$ be a group and let $H \leq G$ have finite index. Then:
\bi
\item[(i)] it is decidable whether or not $H$ is $\V$-closed;
\item[(ii)] $\mathrm{Cl}_{\V}(H)$ is computable from $H$.
\ei
\ep

\proof
(i) Since $[G:H]<\infty$, ${\rm Core}_G(H)$ is a finite intersection of conjugates of $H$ in $G$, and so is computable. In particular $G/{\rm Core}_G(H)$ is computable. Since $\mathbf{V}$ has decidable membership problem, the result follows from Lemma \ref{fi}.

(ii) Since $[G:H] < \infty$, there exist only finitely many subgroups of $G$ containing $H$. Now we decide for each one of them whether it is $\V$-closed or not. The intersection of the $\V$-closed ones must be $\mathrm{Cl}_{\V}(H)$.
\qed

We describe next the $\U$-closed finitely generated subgroups of $F_n$:

\bt
\label{uclosed}
Let $H \leq_{f.g.} F_n$. Then the following conditions are equivalent:
\bi
\item[(i)] $H$ is $\U$-closed;
\item[(ii)] $H$ is $\U$-clopen;
\item[(iii)] $[F_n:H] < \infty$ and $F_n/{\rm Core}_{F_n}(H) \in \U$.
\ei
\et

\proof
(i) $\Rw$ (iii).
Write $F = F_n$. By Proposition \ref{ucl}, if $H$ is $\U$-closed then $\mathrm{Cl}_{\U}(1) \leq H$. We recall that $L_{n,p,p-1} \unlhd F$ for every $p \in \mathbb{P}$. Now Proposition \ref{ucl} yields
$$\mathrm{Cl}_{\U}(1) = \bigcap_{p\in \mathbb{P}} L_{n,p,p-1} \unlhd F.$$

Note that $F^{(2)}$ is a nontrivial normal subgroup of $F$. For every $p \in \mathbb{P}$, we have
$F^{(2)} \leq [K_{n,p-1},K_{n,p-1}] \leq L_{n,p,p-1}$, therefore $F^{(2)} \leq \mathrm{Cl}_{\U}(1) \leq H$.
By \cite[Theorem 1]{KS}, a subgroup of $F$ of infinite index which contains a nontrivial normal subgroup is not finitely generated. 
Hence $[F:H] < \infty$. Finally, $F/\textrm{Core}_F(H) \in \U$ follows from Lemma \ref{fi}.

(iii) $\Rw$ (ii). By  \cite[Proposition 1.2]{MSW}.

(ii) $\Rw$ (i). Trivial.
\qed

\bc
\label{udeccl}
Given $H \leq_{f.g.} F_n$, it is decidable whether or not $H$ is $\U$-closed.
\ec

\proof
Assume that $H$ is given through a finite generating set. As remarked before, it is decidable whether or not $H$ has finite index in $F = F_n$. Assuming now that $[F:H] < \infty$, we can now compute the finitely many conjugates of $H$ \cite[Subsection 2.4]{BS} and compute their intersection \cite[Theorem 2.13]{BS} to get $\textrm{Core}_F(H)$. And we can decide whether or not $F/\textrm{Core}_F(H) \in \U$ by Corollary \ref{mud}. Therefore condition (iii) of Theorem \ref{uclosed} is decidable and the claim follows.
\qed

The next result provides sufficient conditions for a finitely generated subgroup of $F_n$ to have a non finitely generated $\U$-closure.

Fix a basis $\{ a_1,\ldots,a_n\}$ of $F_n$. For $j \in [n]$, unless otherwise stated, we denote by $\pi_j:F_n \to \Z$ the surjective homomorphism defined by $$\pi_j(a_i) = \left\{
\begin{array}{ll}
1&\mbox{ if }i = j\\
0&\mbox{ otherwise.}
\end{array}
\right.$$

\bp
\label{fini}
Let $H \leq_{f.g.} F_n$ be such that $\pi_j(H) = 0$ for some $j \in [n]$. Then:
\bi
\item[(i)]
$[F_n:{\rm Cl}_{\U}(H)] = \infty$
\item[(ii)]
${\rm Cl}_{\U}(H)$ is not finitely generated.
\ei
\ep

\proof
(i) Write $F = F_n$ and $C = {\rm Core}_F(\cl_{\U}(H)) \unlhd F$. Suppose that $[F:{\rm Cl}_{\U}(H)]$ is finite. Then $[F:C]$ is also finite. Write $r = [F:C]$.
Then $a_j^rC = (a_jC)^r = C$ and so $a_j^r \in C$. Choose a prime $p > r$. Then $a_j^r \in C \leq \cl_{\U}(H) \leq HL_{n,p,p-1}$ yields
$$\begin{array}{lll}
r&\in&\pi_j(HL_{n,p,p-1}) = \pi_j(H) + \pi_j(L_{n,p,p-1}) = \pi_j([K_{n,p-1},K_{n,p-1}]K_{n,p-1}^p)\\ &&\\
&=&\pi_j([K_{n,p-1},K_{n,p-1}]) + \pi_j(K_{n,p-1}^p) = \pi_j([F,F]F^{p-1})^p\\ &&\\
&=&p(\pi_j([F,F]) + \pi_j(F^{p-1})) = p(p-1)\Z,
\end{array}$$
contradicting $p > r$. Thus $[F:C] = \infty$.

(ii) By Theorem \ref{uclosed}.
\qed

Next, we consider the related concept of $\U$-denseness:

\bt
\label{udense}
Let $H \leq_{f.g.} F_n$. Then the following conditions are equivalent:
\bi
\item[(i)] $H$ is $\U$-dense;
\item[(ii)] $H$ is ${\bf Su}$-dense.
\ei
\et

\proof
(i) $\Rw$ (ii). By a theorem in \cite{MSTs}, $H$ is ${\bf Su}$-dense if and only if, for every prime $p$, we have 
$HN = F_n$ for every $N \unlhd F_n$ such that $F_n/N \in \U_p$. Since $\U_p \subseteq \U$, the latter property follows from $H$ being $\U$-dense.

(ii) $\Rw$ (i). Since $\U \subseteq {\bf Su}$ by Theorem \ref{charU}, the implication follows from \cite[Corollary 3.2]{MSW}.
\qed

\bc
\label{udecd}
Given $H \leq_{f.g.} F_n$, it is decidable whether or not $H$ is $\U$-dense.
\ec

\proof
It was proved in \cite{MSTs} that it is decidable whether or not $H$ is ${\bf Su}$-dense, therefore the result follows from Theorem \ref{udense}.
\qed

The next example shows that at least one of the inclusions in Proposition \ref{ucl} can be strict.

\be
\label{exin}
There exists some $H \unlhd F_n$ of finite index such that $\mathrm{Cl}_{\M}(H) < \mathrm{Cl}_{\U}(H)$.
\ee

Indeed, by Proposition \ref{UinM} there exists some $G \in \M \setminus \U$. Let $\p:F_n \to G$ be a surjective homomorphism and let $H = \ker(\p)$. By Lemma \ref{fi}, $H$ is $\M$-closed but not $\U$-closed. Therefore 
$\mathrm{Cl}_{\M}(H) = H < \mathrm{Cl}_{\U}(H)$.
\qed

The case of free groups of infinite rank is easily settled in the following result:

\bp
\label{infrank}
Let $F$ be a free group of infinite rank and let $H \leq_{f.g.} F$. Then $H$ is neither $\U$-closed nor $\U$-dense.
\ep

\proof
Let $A$ be a basis of $F$. For every subset $B$ of $A$, there exists a homomorphism $\p:F \to F_B$ such that $\p|_{F_B}$ is the identity, i.e. $F_B$ is a {\em retract} of $F_A$. By \cite[Proposition 1.6]{MSW}, the pro-$\U$ topology on $F_B$ is the subspace topology with respect to the pro-$\U$ topology on $F$. It follows from general topology that
\beq
\label{infrank1}
\cl_{\U}^{F_B}(X) = \cl_{\U}^{F}(X) \cap F_B\mbox{ for all $B \subseteq A$ and }X \subseteq F_B.
\eeq
Since $H$ is finitely generated and $A$ is infinite, there exists some finite $B \subset A$ such that $H \leq F_B$. Let $c \in A \setminus B$ and write $C = B \cup \{ c\}$. 

Suppose that $H$ is $\U$-closed. Then (\ref{infrank1}) yields $\cl_{\U}^{F_C}(H) = H$ and so $[F_C:H] < \infty$ by Theorem \ref{uclosed}. Since $H \leq F_B < F_C$, we get $[F_C:F_B] < \infty$, a contradiction. Therefore $H$ is not $\U$-closed.

We define next a homomorphism $\theta:F \to C_2$ as follows. For every $a \in A$, let
$$\theta(a) = \left\{
\begin{array}{ll}
1&\mbox{ if }a = c\\
0&\mbox{ otherwise.}
\end{array}
\right.$$
Write $K = \ker(\theta)$. Then $K \lhd F$ and $F/K \cong C_2 \in \U$. Since $HK = K < F$, then $H$ is not $\U$-dense.
\qed

A more challenging problem consists of studying the pro-$\U$ topology for finitely generated subgroups of the free object $F_A({\cal{V}}(\U))$ on some set $A$. The following condition (which $\U$ does not satisfy in view of Corollary \ref{novfin}) certainly simplifies the matters at issue:

\bp
\label{provprof}
Let $\V$ be a pseudovariety of finite groups such that $({\cal{V}}(\V))^f = \V$. Then the pro-$\V$ and the profinite topologies on $F_A({\cal{V}}(\V))$ coincide for every set $A$.
\ep

\proof
We know that, for the pro-$\V$ (respectively profinite) topology, $G = F_A({\cal{V}}(\V))$ is a topological group where the normal subgroups $K$ of $G$ such that $G/K\in \mathbf{V}$ (respectively $G/K\in \mathbf{G}$) form a basis of neighbourhoods of the identity. Hence it suffices to show that, for every $K \unlhd G$, we have $G/K\in \mathbf{V}$ if and only if $G/K\in \mathbf{G}$.

The direct implication holds trivially. Conversely, assume that $G/K\in \mathbf{G}$. Since $G \in {\cal{V}}(\V)$, we have $G/K \in {\cal{V}}(\V)$ as well. Now $G/K$ finite yields $G/K \in ({\cal{V}}(\V))^f = \V$ and we are done.
\qed

\bc
\label{DC}
Let $A$ be an arbitrary set. Then:
\bi
\item[(i)] every finitely generated subgroup of $F_A({\cal{V}}(\ab))$ is $\ab$-closed;
\item[(ii)] every finitely generated subgroup of $F_A({\cal{V}}(\M))$ is $\M$-closed.
\ei
\ec

\proof
(i) Since free abelian groups are directs sums of copies of $\mathbb{Z}$, we have 
$$\cal{AB} = {\cal{V}}(\mathbb{Z}) = {\cal{V}}(\ab).$$
In 1998, Delgado proved that free abelian groups of finite rank are LERF (i.e. all their finitely generated subgroups are $\G$-closed) \cite{Del}. The arbitrary rank case follows easily from the finite case \cite[Lemma 4.3]{MSTm}. Now the claim follows from Proposition \ref{provprof}.

(ii) In 2000, Coulbois proved that free metabelian groups of finite rank are LERF \cite{Cou}. The arbitrary rank case follows easily from the finite case \cite[Lemma 4.6]{MSTm}. We remarked before that ${\cal{M}} = {\cal{V}}(F_2({\cal{M}})$. Since $F_2({\cal{M}})$ is LERF, it is in fact {\em residually finite} (i.e. $\{1\}$ is $\G$-closed). Thus, for each $u \in F_2({\cal{M}}) \setminus \{ 1\}$, there exists some $M_u \in \M$ and some homomorphism $\p_u:F_2({\cal{M}}) \to M_u$ such that $\p_u(u) \neq 1$. The homomorphisms $\p_u$ induce a homomorphism 
$$\p:F_2({\cal{M}}) \to \ds \prod_{u\in F_2({\cal{M}}) \setminus \{ 1\}} M_u$$
which is clearly injective, hence $F_2({\cal{M}}) \in {\cal{V}}(\M)$ and so $\cal{M} = {\cal{V}}(\M)$.
Now the claim follows from Proposition \ref{provprof}.
\qed

Although $\ab \subset \U \subset \M$ by Proposition \ref{UinM}, the next example shows that Corollary \ref{DC} cannot be generalized to $\U$:

\be
\label{exu}
There exists a finitely generated subgroup of $F_2({\cal{V}}(\U))$ which is not $\U$-closed.
\ee

Indeed, we know that ${\cal{V}}(\U) = \cal{M}$ by Theorem \ref{vum}. Since the quaternion group $Q_8$ is metabelian and 2-generated, there exists some surjective homomorphism $\p:F_2({\cal{M}}) \to Q_8$. Let $N = \ker(\p)$. Since $N$ is a finite index subgroup of a finitely generated group, it is itself finitely generated. 

Suppose that $N$ is $\U$-closed in $G = F_2({\cal{M}})$. By Lemma \ref{fi}, we get $G/{\rm Core}_{G}(N) \in \U$. Since $N \lhd G$, then $$Q_8 \cong G/N = G/{\rm Core}_{G}(N) \in \U,$$
a contradiction (recall the proof of  Proposition \ref{UinM}). Therefore $N$ is not $\U$-closed. This establishes the claim.

\section*{Acknowledgements}

The authors are grateful to the anonymous referee for his helpful comments.

The first author acknowledges support from the Centre of Mathematics of the University of Porto, which is financed by national funds through the Funda\c c\~ao para a Ci\^encia e a Tecnologia, I.P., under the project with references UIDB/00144/2020 and  UIDP/00144/2020.

The second author acknowledges support from the Centre of Mathematics of the University of Porto, which is financed by national funds through the Funda\c c\~ao para a Ci\^encia e a Tecnologia, I.P., under the project with reference UIDB/00144/2020. 

The third author was supported by the Engineering and Physical Sciences Research Council, grant number 
EP/T017619/1.

\vspace{1cm}

{\sc Claude Marion, Centro de
Matem\'{a}tica, Faculdade de Ci\^{e}ncias, Universidade do
Porto, R. Campo Alegre 687, 4169-007 Porto, Portugal}

{\em E-mail address}: claude.marion@fc.up.pt

\bigskip

{\sc Pedro V. Silva, Centro de
Matem\'{a}tica, Faculdade de Ci\^{e}ncias, Universidade do
Porto, R. Campo Alegre 687, 4169-007 Porto, Portugal}

{\em E-mail address}: pvsilva@fc.up.pt

\bigskip

{\sc Gareth Tracey, Mathematics Institute, University of Warwick, Coventry CV4 7AL, U.K.}

{\em E-mail address}: Gareth.Tracey@warwick.ac.uk

\end{document}